\documentclass[11pt,twoside,a4paper,reqno,fleqn]{amsart}
\usepackage{amssymb,amscd}
\allowdisplaybreaks
\usepackage{times}             

\theoremstyle{plain}
\newtheorem{theorem}{Theorem}[section]
\newtheorem{proposition}[theorem]{Proposition}
\newtheorem{lemma}[theorem]{Lemma}

\theoremstyle{definition}
\newtheorem{definition}[theorem]{Definition}

\theoremstyle{remark}
\newtheorem{corollary}[theorem]{Corollary}
\newtheorem{remark}[theorem]{Remark}

\numberwithin{equation}{section}

\newcommand{\R}{{\mathbb R}}
\newcommand{\Rp}{\R_+}
\newcommand{\oRp}{\overline{\R}_+}
\newcommand{\N}{{\mathbb N}}
\newcommand{\Z}{{\mathbb Z}}
\newcommand{\C}{{\mathbb C}}

\newcommand{\bH}{{\mathbb H}}

\newcommand{\cord}{\operatorname{c-ord}}
\newcommand{\As}{\underline{\operatorname{As}}}
\newcommand{\cA}{{\mathcal A}}
\newcommand{\cC}{{\mathcal C}}

\newcommand{\cE}{{\mathcal E}}
\newcommand{\cH}{{\mathcal H}}
\newcommand{\cI}{{\mathcal I}}
\newcommand{\cJ}{{\mathcal J}}
\newcommand{\cL}{{\mathcal L}}
\newcommand{\cM}{{\mathcal M}}
\newcommand{\cO}{{\mathcal O}}

\newcommand{\bR}{{\mathfrak R}}
\newcommand{\bS}{{\mathfrak S}}
\newcommand{\bT}{{\mathfrak T}}
\newcommand{\bm}{{\mathfrak m}}
\newcommand{\br}{{\mathfrak r}}
\newcommand{\bs}{{\mathfrak s}}
\newcommand{\bt}{{\mathfrak t}}
\newcommand{\as}{\textnormal{as}}

\newcommand{\fin}{{\textnormal{fin}}}
\newcommand{\nor}{{\textnormal{nor}}}

\newcommand{\Diff}{\operatorname{Diff}}
\newcommand{\DiffF}{\operatorname{Diff}_{\textnormal{cone}}}
\newcommand{\Symb}{\operatorname{Symb}}

\newcommand{\loc}{\textnormal{loc}}
\newcommand{\comp}{\textnormal{comp}}
\newcommand{\cl}{{\textnormal{cl}}}

\newcommand{\pp}{\ast}
\renewcommand{\max}{\textnormal{max}}
\renewcommand{\min}{\textnormal{min}}
\newcommand{\spa}{\operatorname{span}}
\newcommand{\ds}{{\dim X/2-\delta}}

\newcommand{\bbC}{\boldsymbol{C}}
\newcommand{\bbI}{\boldsymbol{I}}
\newcommand{\bbJ}{\boldsymbol{J}}
\newcommand{\bbTh}{\boldsymbol{\Theta}}

\newcommand{\Res}{\operatorname{Res}}
\renewcommand{\Re}{\operatorname{Re}}

\begin{document}

\title{Green's formulas for cone differential operators}
\author{Ingo Witt} 
\address{University of Potsdam, Institute of Mathematics, 
  D-14415 Potsdam, Germany}
\email{ingo@math.uni-potsdam.de} 
\date{September 10, 2003}

\begin{abstract}
  Green's formulas for elliptic cone differential operators are established.
  This is done by an accurate description of the maximal domain of an elliptic
  cone differential operator and its formal adjoint, thereby utilizing the
  concept of a discrete asymptotic type. From this description, the singular
  coefficients replacing the boundary traces in classical Green's formulas are
  deduced.
\end{abstract}

\renewcommand{\subjclassname}{\textup{2000} Mathematics Subject
     Classification}
\maketitle

\setcounter{tocdepth}{1}
\tableofcontents


\section{Introduction}


\subsection{The main result}

Let $X$ be a compact $C^\infty$--manifold with non-empty boundary, $\partial
X$. On the interior $X^\circ := X\setminus \partial X$, we consider
differential operators $A$ which on $U\setminus \partial X$ for some collar
neighborhood $U\cong [0,1)\times Y$ of $\partial X$, with coordinates $(t,y)$
and $Y$ being diffeomorphic to $\partial X$, take the form
\begin{equation}\label{opera}
  A =t^{-\mu} \sum_{j=0}^\mu a_j(t,y,D_y) (-t\partial_t)^j,
\end{equation}
where $a_j\in C^\infty([0,1);\Diff^{\mu-j}(Y))$ for $0\leq j\leq \mu$.
Such differential operators $A$ are called cone-degenerate, or being of Fuchs
type; written as $A\in \DiffF^\mu(X)$. They arise, e.g., when polar
coordinates are introduced near a conical point.

Throughout, we shall fix some reference weight $\delta\in\R$. This means that
we will be working in the weighted $L^2$--space $\bH^{0,\delta}(X)$ as basic
function space, cf.~\eqref{sc_p1} and also Appendix~\ref{appB}. Let $A^\ast\in
\DiffF^\mu(X)$ be the formal adjoint to $A$, i.e,
\begin{equation}\label{dual}
  \left(Au,v\right) = \left(u,A^\ast v\right), \quad u,\,v\in 
  C_\comp^\infty(X^\circ),
\end{equation}
where $(\;,\,)$ denotes the scalar product in $\bH^{0,\delta}(X)$. Then it is
customary to consider the maximal and minimal domains of $A$,
\[
  D(A_\max) := \bigl\{u\in \bH^{0,\delta}(X)\bigm|Au\in 
  \bH^{0,\delta}(X)\bigr\},
\]
and $D(A_\min)$ is the closure of $C_\comp^\infty(X^\circ)$ in $D(A_\max)$
with respect to the graph norm. Similarly for $D(A_\max^\ast)$,
$D(A_\min^\ast)$. 

Our basic object of study is the \emph{boundary sesquilinear form}
\begin{equation}\label{bo}
  [u,v]_A := \left(Au,v\right) - \left(u,A^\ast v\right), \quad
  u\in D(A_\max),\, v\in D(A_\max^\ast).
\end{equation}
By virtue of \eqref{dual}, $[u,v]_A=0$ if $u\in D(A_\min)$ or $v\in
D(A_\min^\ast)$. Therefore, the boundary sesquilinear form $[\;,\,]_A$
descents to a sesquilinear form
\begin{equation}\label{green}
  [\;,\,]_A\colon D(A_\max)\big/D(A_\min) \times 
  D(A_\max^\ast)\big/D(A_\min^\ast) \to \C,
\end{equation}
denoted in the same manner. 

Our basic task consists in computing \eqref{green}. The result will be called
a \emph{Green's formula\/} in analogy to the classical situation arising in
mathematical analysis. Assuming ellipticity for $A$, cf. Definition~\ref{ell},
what we will actually do is to compute the value of $[\;,\,]_A$ on
distinguished linear bases of $D(A_\max)\big/D(A_\min)$ and
$D(A_\max^\ast)\big/D(A_\min^\ast)$, respectively.

The starting point is as follows: Assuming ellipticity for $A$, any solution
$u=u(x)$ to the equation $Au=f(x)$ on $X^\circ$ possesses an asymptotic
expansion
\begin{equation}\label{asymp}
  u(x) \sim \sum_p \sum_{k+l=m_p-1} \frac{(-1)^k}{k!}\, 
  t^{-p}\log^k t\,\phi_l^{(p)}(y)\quad \textnormal{as $t\to+0$,}
\end{equation}
where the set $\{p\in\C\,|\,m_p\geq1\}$ is discrete, $\Re p\to -\infty$ as
$|p|\to\infty$ on this set, and $\phi_l^{(p)}\in C^\infty(Y)$ for all $p,\,l$,
provided that the right-hand side $f$ possesses a similar expansion. Introduce
the linear operator $T$ acting on the space of all formal asymptotic
expansions of the form \eqref{asymp} by
\begin{multline}\label{tTt}
  \sum_p \sum_{k+l=\underline{m_p-1}} \frac{(-1)^k}{k!}\,t^{-p}\log^k t\,
  \phi_l^{(p)}(y) \\ \mapsto \sum_p \sum_{k+l=\underline{m_p-2}} 
  \frac{(-1)^k}{k!}\,t^{-p}\log^k t\,\phi_l^{(p)}(y).
\end{multline}
As will be seen,  
\begin{itemize}
\item the quotient $D(A_\max)\big/D(A_\min)$ is finite-dimensional,
\item it consists of finite sums of the form \eqref{asymp}, where $\dim
   X/2-\delta-\mu < \Re p<\dim X/2-\delta$,
\item it is left invariant under the action of $T$. 
\end{itemize}
In particular, $T$ as acting on $D(A_\max)\big/D(A_\min)$ is nilpotent.
Similarly for $D(A_\max^\ast)\big/D(A_\min^\ast)$.

\begin{theorem}\label{main}
Let $A\in\DiffF^\mu(X)$ be elliptic. Then, to each Jordan basis 
\begin{equation}\label{jb1}
  \Phi_1,T\Phi_1,\dots,T^{m_1-1}\Phi_1,\dots,\Phi_e,T\Phi_e,
  \dots,T^{m_e-1}\Phi_e
\end{equation}
of $D(A_\max)\big/D(A_\min)$, there exists a unique Jordan basis
\begin{equation}\label{jb2}
  \Psi_1,T\Psi_1,\dots,T^{m_1-1}\Psi_1,\dots,\Psi_e,T\Psi_e,
  \dots,T^{m_e-1}\Psi_e 
\end{equation}
of $D(A_\max^\ast)\big/D(A_\min^\ast)$ such that, for all $i,\,j,\,r,\,s$,
\begin{equation}\label{g-f}
   [T^r \Phi_i,T^s \Psi_j]_A = \begin{cases}
   (-1)^{s+1} & \textnormal{if $i=j$, $r+s=m_i-1$,} \\
   0 & \textnormal{otherwise.}
   \end{cases}
\end{equation}
\end{theorem}

\begin{corollary}
\textup{(a)} Both spaces $D(A_\max)\big/D(A_\min)$ and
$D(A_\max^\ast)\big/D(A_\min^\ast)$ have the same Jordan structure\/
\textup{(}with respect to $T$\textup{)}.

\textup{(b)} The sesquilinear form $[\;,\,]_A$ in \eqref{green} is
non-degenerate.

\textup{(c)} The operator $T$ is skew-adjoint with respect to $[\;,\,]_A$,
i.e.,
\begin{equation}\label{skew}
  [T\Phi,\Psi]_A + [\Phi,T\Psi]_A = 0
\end{equation}
for all $\Phi\in D(A_\max)\big/D(A_\min)$, $\Psi\in
D(A_\max^\ast)\big/D(A_\min^\ast)$.
\end{corollary}

\begin{remark}
The conjugate Jordan basis $\Psi_1,T\Psi_1,\dots,T^{m_1-1}\Psi_1,\dots,\Psi_e,
T\Psi_e,\dots,$ $T^{m_e-1}\Psi_e$ in \eqref{jb2} can be found, at least in
principal, once one controls the first $\mu$ conormal symbols 
$\sigma_c^{\mu-j}(A)(z)$ for $j=0,1,\dots,\mu-1$ of $A$, 
cf.~\eqref{conormal_symbol}. More precisely, let $\bigl\{\bt^{-\mu-k}(z);\,
k\in\N_0\bigr\}$ be the inverse to the complete conormal symbol
$\bigl\{\sigma_c^{\mu-j}(z);\,j\in\N_0\bigr\}$ of $A$ under the Mellin
translation product, cf.~\eqref{MTP}. In particular, $\bt^{-\mu-k}(z)$ for
$k=0,1,2,\dots$ is a meromorphic function on $\C$ taking values in the space
$L_\cl^{-\mu}(Y)$ of classical pseudodifferential operators of order $-\mu$ on
$Y$. Then the Jordan basis in \eqref{jb2} can be computed from the Jordan
basis in \eqref{jb1} and the principal parts of the Laurent expansions of
$\bt^{-\mu-k}(z+\mu)$ around the poles in the strip $\ds-\mu+k< \Re z < \ds$
for $k=0,1,\dots,\mu-1$. See~Theorem~\ref{main1}.
\end{remark}


\subsection{Description of the content}

In Section~\ref{sec2}, we discuss discrete asymptotic types for conormal
asymptotics of the form \eqref{asymp}. The central notions are properness of
an asymptotic type and complete characteristic bases for proper asymptotic
types. In Section~3, we study complete Mellin symbols that form an algebra
unter the Mellin translation product. Here, the main result due to
\textsc{Liu--Witt}~\cite{LW01} states that the type for the asymptotics
annihilated by an elliptic, holomorphic complete Mellin symbol is proper; thus
linking to cone differential operators, cf.~Theorem~\ref{proper}. Then in
Theorem~\ref{main1}, in Section~\ref{sec4}, we establish a formula for the
principal parts of the Laurent expansions around the poles of the inverses to
holomorphic complete Mellin symbols under the Mellin translation product. This
formula involves a complete characteristic basis and its conjugate complete
characteristic basis, similar to the situation arising in Theorem~\ref{main}.
In fact, Theorem~\ref{main1} is one of the two main technical results of this
paper from which Theorem~\ref{main} is easily deduced. The other one is
Theorem~\ref{main2} in Section~\ref{sec5}, where certain ``bi-orthogonality''
relations between the two complete characteristic bases of Theorem~\ref{main1}
are established. Theorem~\ref{main} is proved in Section~\ref{sec6}. We start
with a formula for the boundary sesquilinear form $[\;,\,]_A$ taken from
\textsc{Gil--Mendoza} \cite{GM03}, cf.~Theorem~\ref{gm}. The proof of
Theorem~\ref{main} then consists in evaluating this formula, where the latter
basically means to ``take the residue'' of the formula from
Theorem~\ref{main2}.

In Section~\ref{sec7}, we discuss two examples showing how one can proceed
from the ``Green's formula''~\eqref{g-f} to genuine Green's formulas in
concrete situations. The two appendices are not mandatory for the main text,
but improve understanding. In Appendix~\ref{appA}, we are concerned with local
asymptotic types, i.e., asymptotic types at some fixed singular exponent $p$
from \eqref{asymp}. Already here, all the ingredients from the main text of
the paper occur in embryonic form. For instance, the forerunner of
Theorem~\ref{main2} is a famous formula due to \textsc{Keldysh} \cite{Kel51},
see~Remark~\ref{nmb}~(b). An analogue of the boundary sesquilinear form
$[\;,\,]_A$ is also provided, see~\eqref{harry} and Proposition~\ref{awq}.  In
Appendix~\ref{appB}, we describe $D(A_\max)$, $D(A_\min)$ as function spaces
with asymptotics. Among others, this gives a concise way of identifying the
quotient $D(A_\max)\big/D(A_\min)$.

Let us mention some related work: Green's formulas have been under
investigation for a long period, see, e.g.,
\textsc{Coddington--Levinson}~\cite{CL55} for O.D.E.  and
\textsc{Lions--Ma\-ge\-nes}~\cite{LM68} for P.D.E. For singular situations,
see, e.g., \textsc{Nazarov--Plamenevskij} \cite{NP94}. Our approach to
cone-de\-gen\-er\-ate differential operators is built upon work of
\textsc{Schulze}~\cite{Sch91, Sch98}. For instance, the fact that the quotient
$D(A_\max)\big/D(A_\min)$ is finite-dimensional and consists of formal
asymptotic expansions of finite length is an easy consequence, see also
\textsc{Lesch}~\cite[Section~1.3]{Les97}. Recently,
\textsc{Gil--Mendoza}~\cite{GM03} received results similar to ours. Without
reaching the final formula \eqref{g-f}, they studied much of the structure of
the boundary sesquilinear form \eqref{green}. In case $A$ is symmetric, they
applied their results to describe the self-adjoint extensions of $A$.
Keldysh's formula was thoroughly discussed in
\textsc{Kozlov--Maz'ya}~\cite[Appendix~A]{KM99}.


\subsection{Notation}\label{nota}

Notation introduced here will be used without further comment.

$\bullet$\enspace Scalar products on $L^2(Y)$ are given by
\begin{equation}\label{nuzi}
  \left(\phi,\psi\right) := \int_Y \phi(y)\overline{\psi(y)}\,d\mu(y), \quad
  \left\langle\phi,\psi\right\rangle := \int_Y \phi(y)\psi(y)\,d\mu(y),
\end{equation}
where $d\mu$ is a fixed positive $C^\infty$--density $d\mu$ on $Y$. For an
operator $B$ on $C^\infty(Y)$, its formal adjoint $B^\ast$ is defined with
respect to the scalar product $\left(\;,\,\right)$, while the transpose $B^t$
is defined with respect to $\left\langle\;,\,\right\rangle$. In particular,
\[
  \overline{B^\ast \phi} = B^t \bar\phi, \quad \phi\in C^\infty(Y). 
\]
For $u,\,v\in\bH^{0,\delta}(X)$, one supported in the collar neighborhood
$U$ of $\partial X$,
\begin{equation}\label{sc_p1}
  (u,v) :=  \int_{(0,1)\times Y} t^{\dim X-2\delta-1}\, 
  u(t,y)\overline{v(t,y)}\,dtd\mu(y),
\end{equation}
cf.~\eqref{dual}, \eqref{bo}. There should be no ambiguity of usage the same
symbol $(\;,\,)$ in the two different situations \eqref{nuzi}, \eqref{sc_p1}.

$\bullet$\enspace Let $J$ be a finite-dimensional linear space and $T$ be a
nilpotent operator acting on $J$. Then $\Phi_1,\dots,\Phi_e$ is called a
\emph{characteristic basis\/} of $J$ (with respect to $T$, where the latter is
often understood from the context) if
\begin{equation}\label{jord}
  \Phi_1,T\Phi_1,\dots,T^{m_1-1}\Phi_1,\dots,\Phi_e,
  T\Phi_e,\dots,T^{m_e-1}\Phi_e 
\end{equation}
for certain integers $m_1,\dots,m_e\geq1$ form a linear basis of $J$. The
matrix of $T$ with respect to such a linear basis has Jordan form. Therefore,
a characteristic basis $\Phi_1,\dots,\Phi_e$ always exists, the integers
$m_1,\dots,m_e$ are uniquely determined (up to permutation) and equal the
sizes of the Jordan blocks, and $e$ is the number of the Jordan blocks. The
tuple $(m_1,\dots,m_e)$ is called the \emph{characteristic\/} of $J$ (or of
the characteristic basis $\Phi_1,\dots,\Phi_e$).

$\bullet$\enspace Let $E$ be either the space $C^\infty(Y)$ (in Section~2) or
a Banach space (in Appendix~A). Then $E^\infty:=\bigcup_{m\in\N} E^m$ denotes
the space of finite sequences in $E$, where we identify $E^m$ as linear
subspace of $E^{m+1}$ through the map $(\phi_0,\dots,\phi_{m-1})\mapsto
(0,\phi_0,\dots, \phi_{m-1})$, i.e., by adding a leading zero. For $\Phi\in
E^\infty$, let $m(\Phi)$ be the least integer $m$ so that $\Phi\in E^m$. The
right shift operator $T$ sending $(\phi_0,\phi_1,\dots,\phi_{m-1})$ to
$(\phi_0,\phi_1,\dots,\phi_{m-2})$ acts on $E^\infty$. In particular,
$T^{m(\Phi)}\Phi=0$, while $T^i\Phi\neq0$ for $0\leq i\leq m(\Phi)-1$. (In
case $E=C^\infty(Y)$, the operator $T$ is directly related to the operator
$T$ in \eqref{tTt}, cf.~Remark~\ref{rem1}.)

$\bullet$\enspace For $E$ as above, $p\in\C$, let $\cM_p(E)$ be the space
of germs $E$--valued meromorphic functions and $\cA_p(E)$ be the space of
germs of $E$--valued holomorphic functions at $z=p$. We identify the quotient
$\cM_p(E)\big/\cA_p(E)$ with the space $E^\infty$ through the map
\[
  \frac{\phi_0}{(z-p)^m} + \frac{\phi_1}{(z-p)^{m-1}} + \dots +
  \frac{\phi_{m-1}}{z-p} \mapsto (\phi_0,\phi_1,\dots,\phi_{m-1}).
\]
Then $T$ corresponds to multiplication by $z-p$. For $F\in\cM_p(E)$, let
$[F(z)]_p^\ast$ denote the principal part of $F(z)$ at $z=p$. For
$\Phi=(\phi_0,\phi_1,\dots,\phi_{m-1})\in E^\infty$, we set
\[
  \Phi[z-p]:= \frac{\phi_0}{(z-p)^m}+\frac{\phi_1}{(z-p)^{m-1}}+
  \dots+\frac{\phi_{m-1}}{z-p}\in \cM_p(E).
\]

$\bullet$\enspace Now let $E=C^\infty(Y)$. For $\phi,\,\psi\in C^\infty(Y)$,
let $\phi\otimes \psi$ be the rank-one operator $C^\infty(Y)\ni h\mapsto
(h,\psi)\phi\in C^\infty(Y)$. More generally, for $F,\,G\in
\cM_p(C^\infty(Y))$, we introduce the meromorphic operator family $F(z)\otimes
G(z)$ by
\[
  F(z)\otimes G(z)h := \bigl(h,G(\bar z)\bigr)F(z), \quad
  h\in C^\infty(Y),
\]
where $\bigl(h,G(\bar z)\bigr)=\bigl\langle h,\bar G(z)\bigr\rangle$ and $\bar
G(z):=\overline{G(\bar z)}$. For $\Phi,\,\Psi\in
[C^\infty(Y)]^\infty$, we further set
\begin{multline}\label{sybe1}
  (\Phi\otimes \Psi)[z-p] \\ := \frac{\phi_0\otimes \psi_0}{(z-p)^m} +  
  \frac{\phi_0\otimes \psi_1 + \phi_1\otimes\psi_0}{(z-p)^{m-1}} + \dots
  \frac{\phi_0\otimes \psi_{m-1}+\dots + \phi_{m-1}\otimes\psi_0}{z-p},
\end{multline}
where $\Phi=(\phi_0,\dots,\phi_{m-1})$, $\Psi=(\psi_0,\dots,\psi_{m-1})$. In
particular,
\begin{equation}\label{sybe2}
  (\Phi\otimes\Psi)[z-p]=(z-p)^m\, \bigl[\Phi[z-p]\otimes 
  \Psi[z-p]\bigr]_p^\pp,
\end{equation}
where $m=\max\{m(\Phi),m(\Psi)\}$. 

$\bullet$\enspace If $E$ is a Banach space, then we use the same notation, but
with $(\;,\,)$ replaced with the dual pairing $\langle\;,\,\rangle$ between
$E$ and $E'$. This is due to the circumstance that in this situation we are
working with the dual instead of the anti-dual to $E$. In particular,
$\phi\otimes \psi$ for $\phi\in E$, $\psi\in E'$ means the rank-one operator
$E\ni h \to \langle h,\psi\rangle \phi\in E$, while \eqref{sybe1},
\eqref{sybe2} are formally unchanged.

$\bullet$\enspace On $[C^\infty(Y)]^\infty$, we introduce three commuting
involutions by
\begin{align*}
  \bbC\Phi &:=(\bar\phi_0,\bar\phi_1,\dots,\bar\phi_{m-2},\bar\phi_{m-1}), \\
  \bbI\Phi &:=((-1)^m\bar\phi_0,(-1)^{m-1}\bar\phi_1,\dots,
     \bar\phi_{m-2},-\bar\phi_{m-1}),\\
  \bbJ\Phi &:=((-1)^m\phi_0,(-1)^{m-1}\phi_1,\dots,\phi_{m-2},-\phi_{m-1}),
\end{align*}
where $\Phi=(\phi_0,\phi_1,\dots,\phi_{m-2},\phi_{m-1})$. Note that $\bbI=\bbC
\bbJ$, $T\bbC=\bbC T$, $\bbI T+T\bbI=0$, and $\bbJ T+T\bbJ=0$. 

$\bullet$\enspace The cut-off function $\omega\in C_\comp^\infty(\oRp)$
satisfies $\omega(t)=1$ if $t\leq 1/2$, $\omega(t)=0$ if $t\geq1$. It is used
to localize into the collar neighborhood $U\cong [0,1)\times Y$ of $\partial
X$.


\section{Discrete asymptotic types}\label{sec2}


The notion of discrete asymptotic type for conormal cone asymptotics goes back
to \textsc{Rempel--Schulze} \cite{RS89} in the one-dimensional and
\textsc{Schulze} \cite{Sch91} in the higher-dimen\-sion\-al case. It allows to
integrate asymptotic information into a functional-analytic setting, cf.~also
Appendix~\ref{appB}. The refinements of this notion presented here are due to
\textsc{Liu--Witt} \cite{LW01}.


\subsection{Preliminaries}

Let $C_\as^{\infty,\delta}(X)$ be the space of all $u\in C^\infty(X^\circ)$
possessing an asymptotic expansion as in \eqref{asymp}, as $x\to\partial X$,
where, additionally, $\Re p<\ds$ if $m_p\geq1$. Moreover, let
$C_\cO^\infty(X)$ be the space of all $u\in C^\infty(X^\circ)$ vanishing to
the infinite order on $\partial X$ (i.e., $\phi_l^{(p)}(y)=0$ for all $p,\,l$
in \eqref{asymp}).

Henceforth, we shall fix a splitting $U\setminus\partial X\cong(0,1)\times Y$,
$x\mapsto(t,y)$ of coordinates near $\partial X$. It turns out, however, that
our constructions are independent of this chosen splitting of coordinates,
cf.~Remarks~\ref{rem1}, \ref{ccc}, and \ref{rem3} and Proposition \ref{explai}.

\begin{definition}\label{nmp}
(a) A discrete subset $V\subset\C$ is called a \emph{carrier of asymptotics\/} 
if $|\Re p|\to\infty$ on $V$ as $|p|\to\infty$. For $\delta\in\R$, we write
$V\in\cC^\delta$ if, in addition, $V\subset\{z\in\C;\,\Re z<\dim
X/2-\delta\}$.

(b) For $V\in\cC^\delta$, we define $\cE_V^\delta(Y)$ to be the space of all
mappings $\Phi\colon \C\to [C^\infty(Y)]^\infty$ satisfying $\{p\in
\C\,|\,\Phi(p)\neq0\}\subseteq V$. In particular, we have
$\cE_V^\delta(Y)=\prod_{p\in V}[C^\infty(Y)]_p^\infty$, where
$[C^\infty(Y)]_p^\infty$ is an isomorphic copy of $[C^\infty(Y)]^\infty$.
Moreover, we set $\mathcal{E}^\delta(Y):=
\bigcup_{V\in\cC^\delta}\cE_V^\delta(Y)$.
\end{definition}

The operations $T$, $\bbC$, $\bbI$, and $\bbJ$ are point-wise defined on
$\cE^\delta(Y)$, cf.~Section~\ref{nota}. For instance, $T\Phi(p) = T(\Phi(p))$
for $\Phi\in \cE^\delta(Y)$, $p\in\C$. We also write $m^p(\Phi)$ instead of
$m(\Phi(p))$. 

We next provide an isomorphism 
\begin{equation}\label{non_canon}
  C_{\as}^{\infty,\delta}(X)\big/C_\cO^\infty(X)\to \cE^\delta(Y)
\end{equation}
that is non-canonical in the sense that it depends on the choice of splitting
of coordinates near $\partial X$: With the vector $\Phi\in\cE_V^\delta(Y)$,
where $\Phi(p)= (\phi_0^{(p)},\dots, \phi_{m_p-1}^{(p)})$ for $p\in V$, we
associate the formal asymptotic expansion
\[
  u(x)\sim \sum_{p\in V}\sum_{k+l=m_p-1}\frac{(-1)^k}{k!}\
    t^{-p}\log^k t\ \phi_l^{(p)}(y) \enspace \text{as $t\to+0$,} 
\]
cf.~\eqref{asymp}. (To see that \eqref{non_canon} is surjective needs to
invoke a Borel-type argument.)

\begin{remark}\label{rem1}
The operator $T$ acting on the quotient $C_{\as}^{\infty,\delta}(X)
\big/C_\cO^\infty(X)$, as introduced in \eqref{tTt}, is well-defined, i.e., it
is independent of the chosen splitting of coordinates near $\partial X$.
Moreover, the isomorphism \eqref{non_canon} intertwines this operator and the
right-shift operator $T$ acting on $\cE^\delta(Y)$.
\end{remark}

We need some further notation: For $\Phi\in\cE^\delta(Y)$, we introduce
\[
  \cord(\Phi) := \dim X/2-\max\{\Re p;\,\Phi(p)\neq0\}
\]
(the ``\emph{conormal order\/}'' of $\Phi$ understood in an $L^2$-sense). Note
that $\cord(\Phi)>\delta$ if $\Phi\in\cE^\delta(Y)$, and
$\cord(T^k\Phi)\to\infty$ as $k\to\infty$.  Note also that, for
$\Phi_i\in\cE^\delta(Y)$, $\alpha_i\in\C$ for $i=1,2,\dots$ satisfying
$\cord(\Phi_i)\to\infty$ as $i\to\infty$, the series
$\sum_{i=1}^\infty\alpha_i\Phi_i$ is explained in $\cE^\delta(Y)$ in a natural
fashion. In particular,
\[
  \sum_{i=1}^\infty \alpha_i\Phi_i=0 \quad \Longleftrightarrow \quad 
  \cord\bigl(\sum_{i=1}^{i_0}\alpha_i\Phi_i\bigr)\to\infty \quad 
  \text{as $i_0\to\infty$}.
\]
Furthermore, $\Phi\in\cE^\delta(Y)$ is called a \emph{special vector\/} if
$\Phi\in \cE_{p-\N_0}^\delta(Y)$ for some $p\in\C$. If $\Phi\neq0$, then $p$
is uniquely determined by $\Phi$ and the additional requirement that
$\Phi(p)\neq0$. This complex number $p$ is denoted by $\gamma(\Phi)$.


\subsection{Definition of discrete asymptotic types}

Discrete asymptotic types are certain linear subspace of the space 
$C_\as^{\infty,\delta}(X)\big/C_\cO^\infty(X)$ of all formal asymptotic
expansions. 

\begin{definition}\label{fko}
A \emph{discrete asymptotic type\/}, $P$, for conormal cone asymptotics as 
$x\to\partial X$, of conormal order at least $\delta$, is a linear subspace
of $C_\as^{\infty,\delta}(X)\big/C_\cO^\infty(X)$ that is represented, in the
given splitting $U\setminus\partial X\cong(0,1)\times Y$, $x\mapsto (t,y)$ of
coordinates near $\partial X$, through the isomorphism \eqref{non_canon} by a
linear subspace $J\subset \cE_V^\delta(Y)$ for some $V\in\cC^\delta$
satisfying the following conditions:

(i) $TJ\subseteq J$. 

(ii) $\dim J^{\delta+j}<\infty$ for all $j\in\N_0$, where $J^{\delta'}:=
J/(J\cap \cE^{\delta'}(Y))$ for $\delta'\geq \delta$.

(iii) There is a sequence $\{p_i;\,1\leq j<e+1\}\subset\C$, where $e\in\N_0\cup
\{\infty\}$, such that $\Re p_i<\ds$ for all $i$, $\Re p_i\to-\infty$ as
$i\to\infty$ if $e=\infty$, $V\subseteq \bigcup_{i=1}^e \{p_i\}-\N_0$, and
\[
  J = \bigoplus_{i=1}^e \left(J\cap \cE_{p_i-\N_0}^\delta(Y)\right). 
\]
The \emph{empty asymptotic type\/}, $\cO$, is represented by the trivial
subspace $\{0\}\subset\cE^\delta(Y)$. The set of all asymptotic types of
conormal order at least $\delta$ is denoted by $\As^\delta(Y)$.
\end{definition}

\begin{remark}\label{ccc}
It can be shown that this notion of discrete asymptotic type is independent of
the splitting of coordinates near $\partial X$. The latter means that changing
coordinates $P\subset C_\as^{\infty,\delta}(X)\big/ C_\cO^\infty(X)$ 
is represented by another linear subspace $J'\subset \cE_{V'}^\delta(Y)$ for
some $V'\in\cC^\delta$ that also satisfies (i)~to~(iii) of
Definition~\ref{fko}.
\end{remark}

Let $P,\,P'\in\As^\delta$ be represented by $J,\,J'\subset \cE^\delta(Y)$,
respectively. For $\delta'\geq \delta$, we say that $P,\,P'$ \emph{coincide up
  to the conormal order\/} $\delta'$ if $J^{\delta'}= J'{}^{\delta'}$ as
subspaces of $\cE^\delta(Y)\big/\cE^{\delta'}(Y)$.  Similarly, for $\delta'>
\delta$, we say that $P,\,P'$ coincide up to the conormal order $\delta'-0$ if
$P,\,P'$ coincide up to the conormal order $\delta'-\varepsilon$ for all
$0<\varepsilon \leq \delta'-\delta$.

It is important to observe that the set $\As^\delta(X)$ of asymptotic types is
partially ordered by inclusion of the representing spaces. This partial order
on $\As^\delta(X)$ will be denoted by $\preccurlyeq$. One of the fundamental
principles in constructing asymptotic types obeying certain prescribed
properties ensues from the following:

\begin{proposition}
$(\As^\delta(X),\preccurlyeq)$ is a lattice with the property that each
non-empty subset \textup{(}resp. each bounded subset\textup{)} possesses a
greatest lower bound \textup{(}resp. a least upper bound\textup{)}.
\end{proposition}

As example, consider $P\in\As^\delta$ and let $\delta'\geq\delta$. Then let
$P^{\delta'}\preccurlyeq P$ denote the largest asymptotic type that coincides
with the empty asymptotic type, $\cO$, up to the conormal order $\delta'$.
Similarly, for $\delta'>\delta$, let $P^{\delta'-0}\preccurlyeq P$ denote the
largest asymptotic type that coincides with the empty asymptotic type up to
the conormal order $\delta'-0$. Of course, in this situation it is easy to
provide representing spaces for $P^{\delta'}$ and $P^{\delta'-0}$,
respectively, but in more involved situations such a task might be not that
simple.

\begin{remark}
There is an abstract concept of introducing asymptotic types if a unital
algebra ${\mathfrak M}$ acting on some linear space ${\mathfrak F}$ ``modulo a
distinguished linear subspace ${\mathfrak F}_0$ in the image'' is given.  In
our case, ${\mathfrak M}=\bigcup_{\mu\in\Z}\Symb_M^\mu(Y)$ is the algebra of
complete Mellin symbols, cf.~Section~\ref{sec3}, ${\mathfrak F}=
C_\as^{\infty,\delta}(X)$, and ${\mathfrak F}_0= C_\cO^\infty(X)$. See
\textsc{Witt}~\cite{Wi2000}.
\end{remark}


\subsection{Proper discrete asymptotic types}

Here we investigate properties of linear subspaces $J\subset \cE_V^\delta(Y)$
satisfying (i) to (iii) of Definition~\ref{fko}.

\begin{proposition}\label{mainy}
Let $J\subset\cE_V^\delta(Y)$ be a linear subspace for some $V\in\cC^\delta$. 
Then there are an $e\in\N_0\cup\{\infty\}$ and a sequence $\{\Phi_i;\,1\leq
i<e+1\}$ of special vectors satisfying $\cord(\Phi_i)\to\infty$ as
$i\to\infty$ if $e=\infty$ such that the vectors $T^k\Phi_i$ for
$1\leq i<e+1$, $k\in\N_0$ span the space $J$ if and only if\/
\textup{(i)}~to\/~\textup{(iii)} of\/ \textup{Definition~\ref{fko}} are
fulfilled.
\end{proposition}

For the rest of this section, assume that $J\subset \cE_V^\delta(Y)$ is a
linear subspace satisfying (i)~to~(iii) of Definition~\ref{fko}.
Let $\Pi_j\colon J\to J^{\delta+j}$ be the canonical surjection. For $j'>j$,
there is a natural surjective map $\Pi_{jj'}\colon J^{\delta+j'} \to
J^{\delta+ j}$ such that $\Pi_{jj''}=\Pi_{jj'} \Pi_{j'j''}$ for $j''>j'>j$ and
\[
   \bigl(J,\Pi_j\bigr) = \varprojlim
   \bigl(J^{\delta+j},\Pi_{jj'}\bigr). 
\]
Note that $T\colon J^{\delta+j} \to J^{\delta+j}$ is nilpotent, where the
operator $T$ is induced by $T\colon J\to J$. Let $(m_1^j,\dots,m_{e_j}^j)$ be
the characteristic of $J^{\delta+j}$, cf.~\eqref{jord} and thereafter.

The sequence $\{\Phi_i;\,1\leq i<e+1\}\subset J$ is said to be a
\emph{characteristic basis\/} of $J$ if there are numbers $m_i\in\N_0\cup
\{\infty\}$, $m_i\geq1$, such that $T^{m_i}\Phi_i=0$ if $m_i<\infty$, while
the sequence $\{T^k\Phi_i;\,1\leq i<e+1,\,0\leq k<m_i\}$ forms a basis of $J$.

\begin{proposition}\label{ghjtz}
Let $J\subset\cE_V^\delta(Y)$ be a linear subspace as above and assume that
$\{\Phi_i;\,1\leq i<e+1\}$ is a characteristic basis of $J$. Then the
following conditions are equivalent:

\textup{(a)} For each $j$, $\{\Pi_j\Phi_1,\dots,\Pi_j\Phi_{e_j}^j\}$ is 
a characteristic basis of $J^{\delta+j}$.

\textup{(b)} For each $j$, $T^{m_1^j-1}\Phi_1,\dots,T^{m_{e_j}-1} \Phi_{e_j}$
are linearly independent over the space $\cE^{\delta+j}(Y)$, while
$T^k\Phi_i\in \cE^{\delta+j}(Y)$, where either $1\leq i\leq e_j$, $k\geq
m_i^j$ or $i>e_j$.

If these conditions hold, then the numbering within the tuples
$(m_1^j,\dots,m_{e_j}^j)$ can be chosen in such a way that, for each $j\geq1$,
there is a characteristic basis $\Phi_1^j,\dots,\Phi_{e_j}^j$ of
$J^{\delta+j}$ such that, for all $j'>j$,
\[
  \Pi_{jj'}\Phi_i^{j'} = \begin{cases}
     \Phi_i^j & \text{if $1\leq i\leq e_j$,} \\
     0 & \text{if $e_j+1\leq i\leq e_{j'}$.}
  \end{cases}
\]
Furthermore, the scheme
\begin{equation} 
\begin{array}{rr} 
\textnormal{$e_1$ rows} &\begin{cases}
\\
\\
\\
\end{cases}\\
\textnormal{$e_2-e_1$ rows} &\begin{cases}
\\
\\
\\
\end{cases}\\
\textnormal{$e_3-e_2$ rows} &\begin{cases}
\\
\\
\\
\end{cases}\\
\\
\end{array}
\begin{matrix}
  m_1^1 & m_1^2 & m_1^3 & m_1^4 & \dots \\
  \hdotsfor{4} & \dots \\
  m_{e_1}^1 & m_{e_1}^2 & m_{e_1}^3 & m_{e_1}^4 & \dots \\[1.5ex]
  & m_{e_1+1}^2 & m_{e_1+1}^3 & m_{e_1+1}^4 & \dots \\
  & \hdotsfor{3} & \dots \\
  & m_{e_2}^2 & m_{e_2}^3 & m_{e_2}^4 & \dots \\[1.5ex]
  & & m_{e_2+1}^3 & m_{e_2+1}^4 & \dots \\
  & & \hdotsfor{2} & \dots \\
  & & m_{e_3}^3 & m_{e_3}^4 & \dots \\[1.5ex]
  & & & \vdots & \ddots, 
\end{matrix}
\end{equation}
where in the $j$th column the characteristic of the space $J^{\delta+j}$
appears, is uniquely determined up to permutation of the $k$th and the $k'$th
row, where $e_j+1\leq k,\,k'\leq e_{j+1}$ for some $j$
\textup{(}$e_0=0$\textup{)}.
\end{proposition}

\begin{definition}
An asymptotic type $P$ is said to be \emph{proper\/} if its representing 
space~$J$ possesses a characteristic basis $\{\Phi_i;\,1\leq i<e+1\bigr\}$
consisting of special vectors that fulfill the equivalent conditions of
Proposition~\ref{ghjtz}. If the tuples $(m_1^j,\dots,m_{e_j}^j)$ are
re-ordered according to this proposition, then the sequence
\begin{equation}
  \bigl\{\bigl(\gamma(\Phi_i);
  m_i^{j_i},m_i^{j_i+1},m_i^{j_i+2},\dots\bigr);\,1\leq i<e+1\bigr\}
\end{equation}
is called the \emph{characteristic\/} of $P$.
\end{definition}

\begin{remark}\label{rem3}
The characteristic 
of a proper asymptotic type $P\in\As^\delta(Y)$ is independent of the
splitting of coordinates near $\partial X$.
\end{remark}

An asymptotic type need not be proper. For an example, see \textsc{Liu--Witt}
\cite[Example~2.23]{LW01}. However, we have the following result, which will be
generalized in Theorem~\ref{proper} below:

\begin{theorem}\label{proper0}
Let $A\in\DiffF^\mu(X)$ be an elliptic cone-degenerate differential operator.
Then
\begin{equation}\label{kill}
  \bigl\{u\in C_\as^{\infty,\delta}(X)\bigm|
  Au \in C_\cO^\infty(X) \bigr\} \big/ C_\cO^\infty(X)
\end{equation}
is a proper asymptotic type. 
\end{theorem}


\section{The algebra of complete Mellin symbols}\label{sec3}


We study the algebra of complete Mellin symbols under the Mellin translation
product. Furthermore, we introduce the important notion of a complete
characteristic basis for the asymptotics annihilated by a holomorphic complete
Mellin symbol.


\subsection{Cone differential operators}

Recall that we have fixed a splitting of coordinates $U\to [0,1)\times Y$,
$x\mapsto (t,y)$ near $\partial X$, $U$ being a collar neighborhood of
$\partial X$. Let $(\tau,\eta)$ be the covariables to $(t,y)$. The compressed
covariable $t\tau$ to $t$ is denoted by $\tilde\tau$, i.e., $(\tilde\tau,
\eta)$ is the linear variable in the fiber of the compressed cotangent bundle
$\tilde T^\ast U$. 

For $A\in \DiffF^\mu(X)$ as given in \eqref{opera}, we denote by
$\sigma_\psi^\mu(A)$ its principal symbol, by $\tilde \sigma_\psi^\mu(A)$ its
compressed principal symbol defined on $\tilde T^\ast U$ by
\[
  \sigma_\psi^\mu(A)(t,y,\tau,\eta) = t^{-\mu}\tilde
  \sigma_\psi^\mu(A)(t,y,t\tau,\eta), \quad 
  (t,y,\tau,\eta)\in T^\ast (U\setminus\partial X)\setminus0,
\]
and by $\sigma_c^{\mu-j}(A)(z)$ for $j=0,1,2,\dots$ its \emph{$j$\/th conormal
  symbol},
\begin{equation}\label{conormal_symbol} 
  \sigma_c^{\mu-j}(A)(z) :=
  \sum_{k=0}^\mu \frac{1}{j!}\,
  \frac{\partial^j a_k}{\partial t^j} (0,y,D_y)z^k, \quad z\in\C.
\end{equation}
$\tilde\sigma_\psi^\mu(A)(t,y,\tilde\tau,\eta)$ is smooth up to $t=0$ and
$\sigma_c^{\mu-j}(z)$ for $j=0,1,2,\dots$ is a holomorphic function in $z$
taking values in $\Diff^\mu(Y)$.

Furthermore, if $A\in\DiffF^\mu(X)$, $B\in\DiffF^\nu(X)$, then
$AB\in\DiffF^{\mu+\nu}(X)$ and
\[
  \sigma_c^{\mu+\nu-l}(AB)(z) = \sum_{j+k=l}\sigma_c^{\mu-j}(A)(z+\nu-k)
  \sigma_c^{\nu-k}(B)(z)
\]
for $l=0,1,2,\dots$ This formula is called the \emph{Mellin translation
  product}.

\begin{definition}\label{ell}
(a) The operator $A\in \DiffF^\mu(X)$ is called \emph{elliptic\/} if $A$ is an 
elliptic differential operator on $X^\circ$ and
\begin{equation}\label{ell1}
  \tilde\sigma_\psi^\mu(A)(t,y,\tilde\tau,\eta) \neq0, \quad
  (t,y,\tilde\tau,\eta)\in \tilde T^\ast U\setminus0. 
\end{equation}

\noindent
(b) The operator $A\in \DiffF^\mu(X)$ is called \emph{elliptic with respect to 
the weight\/} $\delta\in\R$ if $A$ is elliptic in the sense of (a) and, 
in addition,
\begin{equation}\label{fgf}
  \sigma_c^\mu(A)(z)\colon H^{s}(Y)\to H^{s-\mu}(Y), \quad
  \Re z = \ds, 
\end{equation}
is invertible for some $s\in\R$ (and then for all $s\in\R$).
\end{definition}

\begin{proposition}\label{kkg}
If $A\in\DiffF^\mu(X)$ is elliptic, then the set 
\[
  \bigl\{z\in\C\bigm|\text{$\sigma_c^\mu(A)(z)$ regarded as operator in 
  \eqref{fgf} is not invertible}\bigr\}
\]
is a carrier of asymptotics. In particular, there is a discrete set
$D\subset\R$ such that $A$ is elliptic with respect to the weight $\delta$ for
all $\delta\in\R\setminus D$.
\end{proposition}


\subsection{Meromorphic Mellin symbols}

Here we consider the class of meromorphic operator-valued functions arising
in point-wise inverting elliptic conormal symbols $\sigma_c^\mu(A)(z)$. For
further details, see \textsc{Schulze} \cite{Sch98}.

\begin{definition}\label{mellss}
For $\mu\in\Z\cup\{-\infty\}$, the space $\cM_{\as}^\mu(Y)$ of \emph{Mellin 
symbols\/} of order $\mu$ is defined as follows:

(a) The space $\cM_{\cO}^\mu(Y)$ of \emph{holomorphic Mellin symbols\/} of
order $\mu$ is the space of all $L_{\cl}^\mu(Y)$--valued holomorphic functions
$\bm(z)$ on $\C$ such that $\bm(z)\big|_{z=\beta+i\tau}\in
L_{\cl}^\mu(Y;\R_\tau)$ uniformly in $\beta\in[\beta_0,\beta_1]$, where
$-\infty<\beta_0<\beta_1<\infty$.

(b) $\cM_{\as}^{-\infty}(Y)$ is the space of all meromorphic functions
$\bm(z)$ on $\C$ taking values in $L^{-\infty}(Y)$ satisfying the following
conditions:

(i) The Laurent expansion around each pole $z=p$ of
$\bm(z)$ has the form
\begin{equation}\label{fpe}
  \bm(z) = \frac{\bm_0}{(z-p)^\nu}+\frac{\bm_1}{(z-p)^{\nu-1}}+\dots
  + \frac{\bm_{\nu-1}}{z-p} + \sum_{j\geq0} \bm_{\nu+j}(z-p)^j,
\end{equation}
where $\bm_0,\bm_1,\dots,\bm_{\nu-1}\in L^{-\infty}(Y)$ are finite-rank
operators.

(ii) If the poles of $\bm(z)$ are numbered in a certain way,
$p_1,p_2,p_3,\dots$, then $|\Re p_j|\to\infty$ as $j\to\infty$ if the number
of poles is infinite.

(iii) For any function $\chi(z)\in C^\infty(\C)$ such that $\chi(z)=0$ if
$\operatorname{dist}(z,\bigcup_j\{p_j\})\leq1/2$ and $\chi(z)=1$ if
$\operatorname{dist}(z,\bigcup_j\{p_j\})\geq1$, we have
$\chi(z)\bm(z)\big|_{z=\beta+i\tau}\in L^{-\infty}(Y;\R_\tau)$ uniformly in
$\beta\in[\beta_0,\beta_1]$, where $-\infty<\beta_0<\beta_1<\infty$.

(c) We eventually set $\cM_{\as}^\mu(Y):=\cM_\cO^\mu(Y)+
\cM_{\as}^{-\infty}(Y)$. 
\end{definition}

Write $\bm\in\cM_{\as}^\mu(Y)$ for $\mu\in\Z$ as $\bm(z)=\bm_0(z)+ \bm_1(z)$,
where $\bm_0\in \cM_\cO^\mu(Y)$, $\bm_1\in \cM_\as^{-\infty}(Y)$. Then the
(parameter-dependent) principal symbol $\sigma_\psi^\mu\bigl(
\bm_0(z)\big|_{z=\beta+i\tau}\bigr)\in S_\cl^{(\mu)}((T^\ast Y\times\R_\tau)
\setminus 0)$ is independent of the choice of the decomposition of $\bm(z)$
and also independent of $\beta\in\R$.

\begin{definition}\label{elli}
$\bm\in \cM_{\as}^\mu(Y)$ for $\mu\in\Z$ is called \emph{elliptic\/} if
$\sigma_\psi^\mu\bigl(\bm_0(z)\big|_{z=\beta+i\tau}\bigr)\neq0$ 
everywhere.
\end{definition}

\begin{proposition}
\textup{(a)} $\bigcup_{\mu\in\Z}\cM_\as^\mu(Y)$ is a filtered algebra with
respect to the point-wise product as multiplication.

\textup{(b)} $\bm\in \cM_{\as}^\mu(Y)$ is invertible within this algebra,
i.e., with its inverse belonging to $\cM_{\as}^{-\mu}(Y)$, if and only if
$\bm(z)$ is elliptic.
\end{proposition}


We further introduce the algebra $\Symb_M^\mu(Y)$ of \emph{complete Mellin
symbols}. 

\begin{definition}
For $\mu\in\Z\cup\{-\infty\}$, the space $\Symb_M^\mu(Y)$ consists of all
sequences $\mathfrak{S}^\mu=\{\mathfrak{s}^{\mu-j}(z);$ $j\in\N_0\}\subset
\cM_{\as}^\mu(Y)$. Moreover, an element $\mathfrak{S}^\mu\in \Symb_M^\mu(Y)$
is called \emph{holomorphic\/} if $\mathfrak{S}^\mu=\{
\mathfrak{s}^{\mu-j}(z); j\in\N_0\}\subset \cM_\cO^\mu(Y)$.
\end{definition}

\begin{proposition}\label{greif}
\textup{(a)} $\bigcup_{\mu\in\Z}\Symb_M^\mu(Y)$ is a filtered algebra with
involution with respect to the following operations:

\textup{(i)} The Mellin translation product $\bS^\mu\circ_M \bT^\nu
=\{\mathfrak{u}^{\mu+\nu-l}(z);\,l\in\N_0\}\in \Symb_M^{\mu+\nu}$ $(Y)$ for
$\bS^\mu =\{\bs^{\mu-j}(z);\,j\in\N_0\}\in \Symb_M^\mu(Y)$, $\bT^\nu=
\{\mathfrak{t}^{\nu-k}(z);\,k\in\N_0\}\in \Symb_M^\nu(Y)$, where
\begin{equation}\label{MTP}
  \mathfrak{u}^{\mu+\nu-l}(z) = \sum_{j+k=l}\bs^{\mu-j}(z+\nu-k)
  \mathfrak{t}^{\nu-k}(z), \quad l=0,1,2,\dots,
\end{equation}
as multiplication.

\textup{(ii)} The operation $(\bS^{\mu})^{\ast_M}=\{\br^{\mu-j}(z);
\,j\in\N_0\}\in\Symb_M^\mu(Y)$ for $\bS^\mu =\{\bs^{\mu-j}(z);$ $j\in\N_0\}
\in\Symb_M^\mu(Y)$, where
\begin{equation}\label{kily}
  \br^{\mu-j}(z) = \bs^{\mu-j}(\dim X-2\delta-\bar z-\mu+j)^\ast, \quad
  j=0,1,2,\dots,
\end{equation}
as involution.

\textup{(b)} The complete Mellin symbol $\{\bs^{\mu-j}(z);\,j\in\N_0\}\in
\Symb_M^\mu(Y)$ is invertible within the filtered algebra $\bigcup_{\mu\in\Z}
\Symb_M^\mu(Y)$, i.e., with its inverse belonging to $\Symb_M^{-\mu}(Y)$, if
and only if $\bs^\mu(z)$ is elliptic in the sense of\/
\textup{Definition~\ref{elli}}.

\textup{(c)} The map
\begin{equation}\label{mapp}
  \bigcup_{\mu\in\N_0}\DiffF^\mu(X) \to \bigcup_{\mu\in\Z}\Symb_M^\mu(Y), 
  \quad A \mapsto \bigl\{\sigma_c^{\mu-j}(A);\,j\in\N_0\bigr\},
\end{equation}
is a homomorphism of filtered algebras with involution.  
\end{proposition} 


\subsection{The space \boldmath $L_{\bS^\mu}^\delta$}

For a complete Mellin symbol $\bS^\mu$, we introduce a special notation for
the representing space of the ``asymptotics annihilated'' by $\bS^\mu$. We
actually restrict ourselves to holomorphic complete Mellin symbol 
(although the definition can be generalized to meromorphic complete Mellin
symbol by taking into account the possible ``production of asymptotics,''
cf.~\textsc{Liu--Witt} \cite{LW01}). The reason for that is that
Theorem~\ref{proper}, in general, fails to hold without assuming holomorphy.

\begin{definition}
Let $\bS^\mu\in\Symb_M^\mu(Y)$ be holomorphic. Then the linear space 
$L_{\bS^\mu}^\delta\subset\cE^\delta(Y)$ is spanned by all special vectors 
$\Phi\in\cE^\delta(Y)$ satisfying
\begin{equation}\label{cft}
  \sum_{j+k=l} \bs^{\mu-k}(z+k)\Phi(p-j)[z-p+l]\in 
  \cA_{p-l}(C^\infty(Y))
\end{equation}
for $l=0,1,2,\dots$, where $p=\gamma(\Phi)$.
\end{definition}

We also introduce a notation for the expression occuring on the left-hand side
of \eqref{cft}:
\begin{equation}\label{cft2}
  \bbTh_l(\Phi)[z]=
  \bbTh_l(\Phi;{\bS^\mu})[z]:=\sum_{j+k=l} \bs^{\mu-k}(z+k)\Phi(p-j)[z-p+l].
\end{equation}

\begin{proposition}
For $A\in \DiffF^\mu(Y)$, the subspace of $C_\as^{\infty,\delta}(Y)\big/
C_\cO^\infty(X)$ from \eqref{kill} is represented by the linear space
$L_{\bS^\mu}^\delta$, where $\bS^\mu=\bigl\{\sigma_c^{\mu-j}(z);\,
j\in\N_0\bigr\}$.
\end{proposition}

In this situation we sometimes write $L_A^\delta$ instead of
$L_{\bS^\mu}^\delta$.

Generalizing Theorem~\ref{proper0}, we have:

\begin{theorem}[\textsc{Liu--Witt} {\cite[Theorem 2.31]{LW01}}]\label{proper}
For an elliptic, holomorphic $\bS^\mu\in\Symb_M^\mu(Y)$, $L_{\bS^\mu}^\delta$
represents a proper asymptotic type.
\end{theorem}

We also need:

\begin{lemma}
The adjoint relation to \eqref{cft}, 
\[
  \sum_{j+k=l} \br^{\mu-k}(z+k)\Psi(q-j)[z-q+l]\in \cA_{q-l}(C^\infty(Y)),
\]
where $\bR^\mu\in\Symb_M^\mu(Y)$ is as in \eqref{kily}, is equivalent to
\begin{equation}\label{laube}
  \sum_{j+k=l} \bs^{\mu-k}(z)^t\, \bbI\Psi(q-j)[z-p-l]\in 
  \cA_{p+l}(C^\infty(Y)),  
\end{equation}
where $q=\dim X-2\delta-\bar p-\mu$.
\end{lemma}


\subsection{Complete characteristic bases}

The control of asymptotics of the form \eqref{asymp}, of conormal order at
least $\delta$, is equivalent to the control of the conormal symbols
$\sigma_c^{\mu-j}(A)(z)$ for $j=0,1,2,\dots$ in the half-spaces $\Re z<\ds-j$.
We now investigate what happens as $\delta\to-\infty$.
 
Let 
\begin{equation}\label{grew}
  \cE(Y):=  \bigcup_{\delta\in\R}\cE^\delta(Y),
\end{equation}
and $L_{\bS^\mu}:= \bigcup_{\delta\in\R}L_{\bS^\mu}^{\delta}$ for a
holomorphic $\bS^\mu\in\Symb_M^\mu(Y)$.

\begin{definition}
A \emph{complete characteristic basis\/} of $L_{\bS^\mu}$ is the inductive
limit
\[ 
  \varinjlim\left(\{\Phi_h^{\delta};\,h\in\cI^\delta\},
  \tau_{\delta'\delta}\right) 
\] 
of the following inductive system:

(a) For each $\delta\in\R$, $\{\Phi_h^\delta;\,h\in\cI^\delta\}$ is a
characteristic basis of $L_{\bS^\mu}^{\delta}$ of characteristic
$\bigl(\gamma(\Phi_h^{\delta});m_h^{j_h,\delta},$ \linebreak
$m_h^{j_h+1,\delta},\dots\bigr)$ satisfying conditions (a), (b) of
Proposition~\ref{ghjtz}.

(b) For all $\delta>\delta'$, $\tau_{\delta'\delta} \colon \cI^\delta\to 
\cI^{\delta'}$ is an injection such that, for any $h\in\cI^\delta$, 

(i) $\gamma(\Phi_{h'}^{\delta'})=\gamma(\Phi_h^\delta) +a$ for some
$a\in\N_0$,

(ii) $\Phi_h^\delta = T^{m_{h'}^{j_{h'}+a-1,\delta'}} \Phi_{h'}^{\delta'}$,

where $h'=\tau_{\delta'\delta}(h)$ (and $\tau_{\delta''\delta}=
\tau_{\delta''\delta'}\tau_{\delta'\delta}$ for $\delta>\delta'>\delta''$.)
\end{definition}

We write $\cI:= \varinjlim\left(\cI^{\delta},\tau_{\delta'\delta}\right)$ with
injections $\tau_\delta\colon \cI^\delta\to\cI$ (such that $\tau_\delta=
\tau_{\delta'}\tau_{\delta'\delta}$ for $\delta>\delta'$) and
\[ 
  \{\Phi_h;\,h\in\cI\}=
  \varinjlim\left(\{\Phi_h^\delta;\,h\in\cI^\delta\},
  \tau_{\delta'\delta}\right), 
\] 
where each $\Phi_h$ for $h\in\cI$ is the collection
$\left\{\Phi_{h^\delta}^\delta\right\}$ with $h=\left\{h^\delta\right\}$.

The proof of (a) in the next proposition relies on the finite-dimensionality of
the spaces $J^{\delta+j}$:

\begin{proposition}
\textup{(a)} For each holomorphic $\bS^\mu\in\Symb_M^\mu(Y)$, $L_{\bS^\mu}$
possesses a complete characteristic basis.

\textup{(b)} For any complete characteristic basis $\{\Phi_h;\,h\in\cI\}$ of
$L_{\bS^\mu}$, the expression
\[
  T^{m^{p+l}(\Phi_h)}\Phi_h(p) := T^{m^{p+l}(\Phi_{h^\delta}^\delta)}
  \Phi_{h^\delta}^\delta(p), 
\]
where $p\in\C$, $l\in\N_0$, $h^\delta=\tau_\delta(h)$, and $\delta$ is 
chosen in such a way that 
\[
  \delta\leq \max\bigl\{\min\bigl\{\delta'\in\R\bigm|
  \gamma(\Phi_{h^{\delta'}}^{\delta'})-\Re p\in\N_0\bigr\},
  \dim X/2-(\Re p+l)\bigr\},
\]
is well-defined. \textup{(}If there is no $\delta$ such that $\gamma(
\Phi_{h^{\delta}}^\delta)-\Re p\in\N_0$, then we set
$T^{m^{p+l}(\Phi_h)}\Phi_h(p):=0$.\textup{)}
\end{proposition}


\section{Singularity structure of inverses}\label{sec4}


In the sequel, let $\bS^\mu =\{\bs^{\mu-j}(z);\,j\in\N_0\}\in \Symb_c^\mu(Y)$ 
be holomorphic and elliptic. The inverse to $\bS^\mu$ will then be denoted by 
$\bT^{-\mu}=\{\mathfrak{t}^{-\mu-k}(z);\,k\in\N_0\}\in \Symb_M^{-\mu}(Y)$, 
cf.~Proposition~\ref{greif}~(b). In particular, 
\[
  \sum_{j+k=l} \bt^{-\mu-j}(z+\mu)\bs^{\mu-k}(z+k) = 
  \delta_{0l}\operatorname{id}, \quad l=0,1,2,\dots
\]

Before stating Theorems~\ref{main1} and~\ref{main2}, we simplify the situation 
to be considered in the proofs. Due to the facts that 
\emph{
\begin{itemize}
\item in the process of inverting\/ $\bS^\mu$ with respect to the Mellin
  translation product, the ``production of singularities'' of\/
  $\bt^{-\mu-j}(z+\mu)$ at $z=p$ and\/ $\bt^{-\mu-j'}(z+\mu)$ at $z=p'$,
  respectively, influences each other only if $p-p'\in\Z$,
\item control on the singularity structure of\/ $\bt^{-\mu-j}(z+\mu)$ in the
  half spaces\/ $\Re z<\dim X/2-\delta$ for each\/ $\delta\in\R$ provides
  control on the singularity structure of\/ $\bt^{-\mu-j}(z+\mu)$ in the whole 
  of\/ $\C$,
\end{itemize}}
\noindent
we are allowed to assume the following \emph{model situation\/}: The complete
characteristic basis of $L_{\bS^\mu}$ consists of special vectors $\Phi_i$ for
$1\leq i<e+1$, where $e\in\N_0\cup\{\infty\}$, $\gamma(\Phi_1)=p$, and
\begin{equation}\label{assump}
  \gamma(\Phi_i) = p-l, \quad e_l+1\leq i\leq e_{l+1}
\end{equation}
Then $0=e_0<e_1\leq e_2\leq\dots$ and $e=\max\{e_l;\,l\in\N_0\}$. When
referring to this model situation, we denote
\[
  m_i^{l+1} := m^{p-l}(\Phi_i).
\]

\begin{theorem}\label{main1}
Let $\bS^\mu\in \Symb_c^\mu(Y)$ be holomorphic, elliptic. Then, for each 
complete characteristic basis $\{\Phi_h;\,h\in\cI\}$ of $L_{\bS^\mu}$, there 
is a unique complete characteristic basis $\{\Psi_{h^\ast};\,h^\ast \in
\cI^\ast\}$ of $L_{\bR^\mu}$, where $\bR^\mu$ is given by \eqref{kily}, and a 
bijection $\tau^\ast\colon\cI\to \cI^\ast$ such that, for all $p\in\C$ and
$j=0,1,2,\dots$,
\begin{equation}\label{fundamental}
  \left[\bt^{-\mu-j}(z+\mu)\right]_p^\pp = \sum_h T^{m_h^{p+1}}\Phi_h 
  (p-j) \otimes T^{m_{h^\ast}^{q+j+1}}\bbJ\Psi_{h^\ast}(q)[z-p],
\end{equation}
where $q=\dim X-2\delta-\bar p-\mu$, $h^\ast=\tau^\ast(h)$, $m_h^{p+1}=
m^{p+1}(\Phi_h)$, and $m_{h^\ast}^{q+j+1}= m^{q+j+1}(\Psi_{h^\ast})$.
\end{theorem}
\begin{proof}
We assume the model situation \eqref{assump}.

\textbf{Step 1.}\enspace The elements $\Phi$ of $L_{\bS^\mu}$
($=L_{\bS^\mu}^\delta$) are given by the relations
\begin{equation}\label{darst}
  \Phi(p-l)[z-p+l] = \sum_{j=0}^l \left[\bt^{-\mu-j}(z+\mu)
  \phi^{(p-l+j)}(z)\right]_{p-l+j}^\pp, 
\end{equation}
for $l=0,1,2,\dots$, where $\phi^{(p-j)}=\sum_{r=0}^\infty
\phi_r^{(p-j)}(z-p+j)^r\in \cA_{p-j}(C^\infty(Y))$ for $j=0,1,2,\dots$ See
\textsc{Liu--Witt} \cite{LW01}. In fact, the Taylor coefficients
$\phi_r^{(p-j)}\in C^\infty(Y)$ can be chosen arbitrarily, since only a finite
number of them enters the computation of $\Phi(p-l)$.

Since $\Phi(p-l)\in \spa\bigl\{ \Phi_i(p-l);\,1\leq i\leq e_{l+1}\bigr\}$ for
$\Phi\in L_{\bS^\mu}$, we conclude that
\begin{equation}\label{dadar}
  \left[\bt^{-\mu-j}(z+\mu)\right]_{p-l+j}^\pp = \sum_{i=1}^{e_{l+1}}
  \Phi_i(p-l)\otimes H_i^{(jl)}[z-p+l-j]
\end{equation}
for certain $H_i^{(jl)}=\left(h_{0i}^{(jl)},h_{1i}^{(jl)},\dots,
  h_{m_i^{l+1}-1,i}^{(jl)}\right)\in [C^\infty(Y)]^\infty$, which \emph{a
  priori\/} are of length $m_i^{l+1}$. Employing \eqref{dadar}, we rewrite
\eqref{darst} as
\begin{equation}\label{grap}
  \Phi(p-l) = \sum_{i=1}^{e_{l+1}}\sum_{r=0}^{m_i^{l+1}-1} \left( \sum_{j=0}^l
  \sum_{s=0}^r \bigl(\phi_{r-s}^{(p-l+j)}, h_{si}^{(jl)}\bigr) \right)
  T^r\Phi_i(p-l)
\end{equation}
for $l=0,1,2,\dots$

\textbf{Step 2.}\enspace 
We are going to show that \eqref{grap} provides the unique representation of
$\Phi$ as linear combination of the vectors $T^r\Phi_i$.  More precisely, by
induction on $l=0,1,2,\dots$, we construct functions $h_{ri}^{(l)}\in
C^\infty(Y)$ for $r\geq m_i^l$ such that, for each $l$, 
\begin{equation}\label{grap1}
  \text{$h_{m_i^l,i}^{(l)}$ for all $i$ satisfying $m_i^l<m_i^{l+1}$ 
    are linearly independent}
\end{equation}
and
\begin{equation}\label{grap2}
  h_{ri}^{(jl)} = \begin{cases}
  h_{ri}^{(l-j)} & \text{if $r\geq m_i^{l-j}$,} \\
  0 & \text{otherwise.} 
  \end{cases}
\end{equation}
This means that the coefficient in front of $T^r\Phi_i$ equals $\sum_{j=0}^l
\sum_{s=m_i^j}^r \bigl(\phi_{r-s}^{(p-j)},h_{si}^{(j)}\bigr)$ provided that
$m_i^l\leq r\leq m_i^{l+1}-1$.

After the $l$\/th step, $h_{ri}^{(l')}$ will have been constructed for all
$l',\,i,\,r$ satisfying $l'\leq l$, $1\leq i\leq e_{l+1}$, $m_i^{l'}\leq r\leq
m_i^{l+1}-1$. Moreover, \eqref{grap2} will have been proved for all $i,\,r$
satisfying $1\leq i\leq e_{l+1}$, $0\leq r\leq m_i^{l+1}-1$.

\textbf{Base of induction \boldmath $l=0$:}\enspace 
We set $h_{ri}^{(0)}:=h_{ri}^{(00)}$ for $1\leq i\leq e_1$, $0\leq r\leq
m_i^1-1$.

\textbf{Induction step \boldmath $l'<l\to l$:}\enspace 
We write \eqref{grap} as
\begin{multline*}
  \Phi(p-l) = \sum_{i=1}^{e_l}\sum_{r=0}^{m_i^l-1} 
  \left( \sum_{j=0}^l
  \sum_{s=0}^r \bigl(\phi_{r-s}^{(p-l+j)}, h_{si}^{(jl)}\bigr) 
  \right)T^r\Phi_i(p-l) \\
 +\, \sum_{i=1}^{e_{l+1}}\sum_{r=m_i^l}^{m_i^{l+1}-1} 
  \left( \sum_{j=0}^l
  \sum_{s=0}^r \bigl(\phi_{r-s}^{(p-l+j)}, h_{si}^{(jl)}\bigr) \right)
  T^r\Phi_i(p-l).
\end{multline*}
$\Phi$ regarded as a vector in $L_{\bS^\mu}^\delta$ modulo $L_{\bS^\mu}^{
  \delta'}$ for some $\delta'$ satisfying $\dim X/2-\delta'<p-l< \dim
X/2-\delta'+1$ is a unique linear combination of the vectors $T^r\Phi_i$ for
$1\leq i\leq e_{l+1}$, $0\leq r\leq m_i^{l+1}-1$. The coefficient in front of
$T^r\Phi_i$ for $1\leq i\leq e_l$, $0\leq r\leq m_i^l-1$ is known if one knows
$\Phi$ modulo $L_{\bS^\mu}^{\delta'-1}$. By inductive hypothesis, this
coefficient equals $\sum_{j=0}^{l-1} \sum_{s=m_i^j}^r
\bigl(\phi_{r-s}^{(p-j)},h_{si}^{(j)}\bigr)$. Thus, we obtain \eqref{grap2}
for $1\leq i\leq e_l$, $0\leq r\leq m_i^l-1$, since the functions
$\phi_r^{(p-j)}\in C^\infty(Y)$ are arbitrary.

It remains to set $h_{ri}^{(l')}:=h_{ri}^{(l-l',l)}$ for $0\leq l'\leq l$,
$m_i^l\leq r\leq m_i^{l+1}-1$.

In particular,
\[
  H_i^{(jl)} = \left(h_{m_i^{l-j},i}^{(l-j)},h_{m_i^{l-j}+1i}^{(l-j)},\dots,
  h_{m_i^{l+1}-1,i}^{(l-j)}\right)
\]
is actually of length $m_i^{l+1}-m_i^{l-j}$.

\textbf{Step 3.}\enspace 
We now fix some $a\in\N_0$ and set 
\[
  H_i^{(l)} = H_{i;a}^{(l)}:= \left(h_{m_i^l,i}^{(l)},h_{m_i^{l-j}+1i}^{(l)},
  \dots,h_{m_i^a-1,i}^{(l)}\right)
\]
for $0\leq l\leq a-1$. Then $H_i^{(jl)} = T^{m_i^a-m_i^{l+1}} H_i^{(l-j)}$ and
\begin{multline}\label{gqe}
  \left[\bt^{-\mu-j}(z+\mu-k)\right]_{p-b}^\pp \\ 
  = \,\sum_{i=1}^{e_{l+b+1}} T^{m_i^{k+b}}\Phi_i(p-l-b)\otimes 
  T^{m_i^a-m_i^{l+b+1}} H_i^{(k+b)}[z-p+b]
\end{multline}
for all $j,\,k,\,l,\,b$ satisfying $j+k=l$, $l+b<a$. We shall employ
\eqref{gqe} to show that the vectors $\Psi_1,\dots,\Psi_{e_a}$ defined by
\begin{equation}\label{faqr}
  \Psi_i(q+l) := \bbJ H_i^{(l)}, \quad 0\leq l\leq a-1,
\end{equation}
where $q=\dim X-2\delta-\bar p-\mu$, form a characteristic basis of
$L_{\bR^\mu}^{\delta+\mu-a}$ modulo $L_{\bR^\mu}^{\delta+\mu}$. 

In view of \eqref{grap1}, $\Psi_1,\dots,\Psi_{e_a}$ form a characteristic
basis of the $T$--invariant subspace of $\cE^\delta(Y)$ modulo
$\cE^{\delta+a}(Y)$ generated by these vectors, of characteristic
\[
  \bigl\{\bigl(q+l;m_i^a-m_i^l,m_i^a-m_i^{l-1},\dots,m_i^a-m_i^1,m_i^a
  \bigr);\,1\leq i\leq e_a\bigr\},
\]
where, for a given $i$, $l$ is the least integer such that $m_i^{l+1}=m_i^a$.
In particular, the dimension of this space equals $\dim
L_{\bS^\mu}^\delta\big/L_{\bS^\mu}^{\delta+a}=\sum_{i=1}^{e_a}m_i^a$.
Invoking a duality argument, we see that it suffices to prove that each
$\Psi_i$ belongs to $L_{\bR^\mu}^{\delta+\mu-a}$ modulo
$L_{\bR^\mu}^{\delta+\mu}$.

\textbf{Step 4.}\enspace
By virtue of \eqref{laube}, we have to show that
\begin{equation}\label{hajo}
  \sum_{j+k=l} \bs^{\mu-k}(z)^t\, \bbC H_i^{(a-j-1)}[z-\tilde p-l]
   = O(1) \quad \text{as $z\to \tilde p+l$}
\end{equation}
for $1\leq i\leq e_a$, $0\leq l\leq a-1$, where $\tilde p:=p-a+1$.

For $l+b=a-1$, we infer from \eqref{gqe}
\begin{multline*}
  \delta_{0l} \operatorname{id} = 
  \sum_{j+k=l} \bt^{-\mu-j}(z+\mu-k)\bs^{\mu-k}(z) \\
  = \sum_{j+k=l} \sum_{i=1}^{e_a} (z-\tilde p-l)^{m_i^a-m_i^{a-j-1}}
  \left(T^{m_i^{a-j-1}}\Phi_i(\tilde p)[z-\tilde p-l] \right) \\ 
  \otimes \left(\bs^{\mu-k}(z)^\ast\,
  H_i^{(a-j-1)}[z-\tilde p-l]\right) + O(1) \quad \text{as $z\to \tilde p+l$.}
\end{multline*}
Since the leading entries of the vectors $T^{m_i^{a-j-1}}\Phi_i(\tilde p)$ (if
there are any) for different~$i$ are linearly independent, we arrive at
\eqref{hajo}.

\textbf{Step 5.}\enspace
Returning to the notation $H_i^{(l)}=H_{i;a}^{(l)}$, we see that the $\Psi_i$
defined by \eqref{faqr} for $1\leq i<e+1$ as $a\to\infty$ constitute a
complete characteristic basis of $L_{\bR^\mu}$ modulo
$L_{\bR^\mu}^{\delta+\mu}$. Furthermore, the considerations also show
uniqueness for the complete characteristic basis of $L_{\bR^\mu}$ modulo
$L_{\bR^\mu}^{\delta+\mu}$ just constructed.
\end{proof}


\section{Generalization of Keldysh's formula}\label{sec5}


Conjugacy of complete characteristic bases in the sense of Theorem~\ref{main1}
forces certain \emph{bilinear relations\/} between the bases elements to hold,
as for local asymptotic types. We are now going to derive these relations 
keeping the notations of the previous section.

\begin{theorem}\label{main2}
For all $p,\,h,\,h^\ast,\,l,\,j$ satisfying $j\leq l$,
\begin{multline}\label{graff}
  \sum_{r=j}^l \left\langle\bbTh_{l-r}(T^{m_h^{p+l+1}}\Phi_h)
    [z+r],T^{m_{h^\ast}^{q+1}}\bbI\Psi_{h^\ast}(q-r)[z-p]
  \right\rangle \\*
  = \delta_{hh^\ast}(z-p)^{-(m_h^p-m_h^{p+l+1})} + 
  O((z-p)^{-(m_h^p-m_h^{p+j})}) \quad \textnormal{as $z\to p$,}
\end{multline}
where $q=\dim X-2\delta-\bar p-\mu$, $m_h^{p+j}=m^{p+j}(\Phi_h)$,
$m_{h^\ast}^{q+1}= m^{q+1}(\Psi_{h^\ast})$, 
\[
  \delta_{hh^\ast} = \begin{cases}
  1 & \textnormal{if $h^\ast=\tau^\ast(h)$,} \\
  0 & \textnormal{otherwise,}  
  \end{cases}
\]
and $\bbTh_{l-r}(T^{m_h^{p+l+1}}\Phi_h) [z+r]$ is defined in \eqref{cft2}.
\end{theorem}

\begin{remark}\label{nmb}
\textup{(a)} \eqref{graff} constitutes an asymptotic expansion formula, with
$j=l$ being the basic case and further correction terms added as $j$ is
getting smaller. In Section~\ref{sec7}, we will be in need of the most refined
case $j=0$.

\textup{(b)} In case $l=0$, we recover Keldysh's formula
\begin{multline*}
  \left\langle\bs^\mu(z)\,T^{m_h^{p+1}}\Phi_h(p)[z-p],
  T^{m_{h^\ast}^{q+1}}\bbI\Psi_{h^\ast}(q)[z-p]\right\rangle \\
  = \delta_{hh^\ast}(z-p)^{-(m_h^p-m_h^{p+1})} 
  + O(1) \quad \textnormal{as $z\to p$.}
\end{multline*}
See \eqref{perw}.
\end{remark}

To prove Theorem~\ref{main2} we need the following simple result: 

\begin{lemma}\label{eay}
Assume the model situation \eqref{assump}. Let $a_{ij}(z)\in\cA_{p-l}(\C)$ for 
some $l\in\N_0$. Then
\[
  \sum_{i=1}^{e_l+1}\sum_{j=0}^l a_{ij}(z)\, T^{m_i^j}\Phi_i(p-l)[z-p+l]\in 
  \cA_{p-l}(C^\infty(Y))
\]
if and only if 
\[
  \sum_{r=0}^j (z-p+l)^{m_i^r} a_{ir}(z) = O((z-p+l)^{m_i^{j+1}})
\]
for all $1\leq i\leq e_l+1$, $0\leq j\leq l$.
\end{lemma}

\begin{proof}[Proof of Theorem~\ref{main2}]
We again assume the model situation \eqref{assump}. 

We reenter the scene at formulas \eqref{gqe}. 
Using these formulas, we write
\begin{multline*}
  \bt^{-\mu-j}(z+\mu) \\ = 
  \sum_{b=0}^{a-j-1}\sum_{i=1}^{e_{j+b+1}}
  (T^{m_i^b}\Phi_i(p-j-b))\otimes (T^{m_i^a-m_i^{j+b+1}} H_i^{(b)})[z-p+b]
  + G_j(z)
\end{multline*}
for $0\leq j\leq a-1$, where $G_j(z)$ is holomorphic on the strip $\dim
X/2-\delta-a+j<\Re z<\dim X/2-\delta$. For any $0\leq l'\leq a-1$, we get
\begin{align}
  \delta_{0l'}\operatorname{id} 
  & = \sum_{j+k=l'} \bt^{-\mu-j}(z+\mu+j)\bs^{\mu-k}(z+l') \label{opi} \\
  & = \sum_{j+k=l'} \sum_{b=0}^{a-l'-1}\sum_{i=1}^{e_{l'+b+1}}
  (z-p+l'+b)^{m_i^{l'+b+1}-m_i^{k+b}} \notag \\
  & \qquad\qquad \left(T^{m_i^{k+b}}\Phi_i(p-l'-b)[z-p+l'+b]\right) \notag \\
  & \qquad\qquad\qquad \otimes \left(\bs^{\mu-k}(z+l')^\ast 
      T^{m_i^a-m_i^{l'+b+1}} H_i^{(k+b)}[z-p+l'+b]\right) \notag \\ 
  & \qquad + \sum_{j+k=l'} G_j(z+j)\bs^{\mu-k}(z+l') \notag
\end{align}
We now apply the operator \eqref{opi} to $T^{m_{i'}^{\tilde
    b-l}}\Phi_{i'}(p-\tilde b+l')[z-p+\tilde b]$, where $l\leq \tilde b\leq
a-1$, and then sum up for $l'$ from $0$ to $l$. Since 
\begin{multline*}
  \sum_{j+k\leq l} G_j(z+j)\bs^{\mu-k}(z+j+k)
  T^{m_{i'}^{\tilde b-l}}\Phi_{i'}(p-\tilde b+j+k)[z-p+\tilde b] \\
\begin{aligned}
  &= \sum_{j=0}^l G_j(z+j)\left(\sum_{k=0}^{l-j}\bs^{\mu-k}(z+j+k)
  T^{m_{i'}^{\tilde b-l}}\Phi_{i'}(p-\tilde b+j+k)[z-p+\tilde b]\right) \\
  &= \sum_{j=0}^l G_j(z+j)
  \bbTh_{l-j}(T^{m_{i'}^{\tilde b-l}}\Phi_{i'})[z+j]\in 
  \cA_{p-\tilde b}(C^\infty(Y)). 
\end{aligned}
\end{multline*}
taking the principal value at $z=p-\tilde b$ on both sides of the resulting
equation, we obtain
\begin{multline*}
  T^{m_{i'}^{\tilde b-l}}\Phi_{i'}(p-\tilde b)[z-p+\tilde b] 
  = \sum_{j+k\leq l}\sum_{i=1}^{e_{\tilde b+1}}
  (z-p+\tilde b)^{m_i^{\tilde b+1}-m_i^{k+b}} \\
  \left\langle\bs^{\mu-k}(z+j+k)\,T^{m_{i'}^{\tilde b-l}}\Phi_{i'}(p-b)
  [z-p+\tilde b], 
  T^{m_i^a-m_i^{\tilde b+1}} \bbC H_i^{(k+b)}[z-p+\tilde b]\right\rangle \\
  T^{m_i^{k+b}}\Phi_i(p-\tilde b)[z-p+\tilde b] + O(1) \quad 
  \text{as $z\to p-\tilde b$,}
\end{multline*}
where $\tilde b=b+j+k$. We get
\begin{multline*}
  T^{m_{i'}^{\tilde b-l}}\Phi_{i'}(p-\tilde b)[z-p+\tilde b]
  = \sum_{i=1}^{e_{\tilde b}+1}\sum_{j=0}^l a_{ij}(z)\,
  T^{m_i^{\tilde b-j}}\Phi_i(p-\tilde b)[z-p+\tilde b] + O(1) \\
  = \sum_{i=1}^{e_{\tilde b}+1}\sum_{j=\tilde b-l}^{\tilde b} 
  a_{i,\tilde b-j}(z)\,T^{m_i^j}\Phi_i(p-\tilde b)[z-p+\tilde b] 
  + O(1) \quad \text{as $z\to p-\tilde b$,}
\end{multline*}
where
\begin{multline*}
  a_{ij}(z)= (z-p+\tilde b)^{m_i^{\tilde b+1}-m_i^{\tilde b-j}} \\
  \times \left\langle \bbTh_{l-j}(T^{m_{i'}^{\tilde b-l}} \Phi_{i'})[z+j],
  T^{m_i^a-m_i^{\tilde b+1}} \bbC H_i^{(\tilde b-j)}[z-p+\tilde b]
  \right\rangle.
\end{multline*}
By virtue of Lemma~\ref{eay}, we conclude that, for all $0\leq j\leq l$,
\begin{multline*}
  \sum_{r=j}^l (z-p+l)^{m_i^{\tilde b-r}} a_{ir}(z) \\
  = \delta_{ii'}(z-p+l)^{m_i^{\tilde b-l}}+
  O((z-p+l)^{m_i^{\tilde b-j+1}}) \quad \text{as $z\to p-l$,}
\end{multline*}
i.e.,
\begin{multline*}
  \sum_{r=j}^l \left\langle \bbTh_{l-r}(T^{m_{i'}^{\tilde b-l}}
  \Phi_{i'})[z+r], T^{m_i^a-m_i^{\tilde b+1}} \bbC H_i^{(\tilde b-r)}
  [z-p+\tilde b]\right\rangle \\
  = \delta_{ii'}(z-p+l)^{-(m_i^{\tilde b+1}-m_i^{\tilde b-l})}
  +O((z-p+l)^{-(m_i^{\tilde b+1}-m_i^{\tilde b-j+1})})
  \quad \text{as $z\to p-l$.}
\end{multline*}
In view of \eqref{faqr}, the latter is \eqref{graff} in the model situation
\eqref{assump}.
\end{proof}


\section{The boundary sesquilinear form}\label{sec6}


In this section, we shall prove Theorem~\ref{main}. From \textsc{Gil--Mendoza} 
\cite[Theorem~7.11]{GM03}, we first quote:

\begin{theorem}\label{gm}
For all $u\in D(A_\max)$, $v\in D(A_\max^\ast)$, 
\begin{multline}\label{gm22}
  [u,v]_A = -\sum_{j=0}^{\mu-1}\; \sum_{\dim X/2-\delta-\mu+j 
  <  \Re p < \dim X/2-\delta} \\ \Res_{z=p} \bigl(\sigma_c^{\mu-j}(A)(z)
  (\omega u)\,\tilde{}\,(z), (\omega v)\,\tilde{}\,(\dim X-2\delta
  -\bar z-\mu+j)\bigr),
\end{multline}
where $\omega(t)$ is a cut-off function and $(\omega u)\,\tilde{}\,(z,\cdot)= 
M_{t\to z}\{(\omega u)(t,\cdot)\}$ denotes the Mellin transform, see around
\eqref{mttt}.
\end{theorem}

\begin{proof}[Proof of Theorem~\ref{main}]
We divide the proof into several steps.

\textbf{Step 1.}\enspace 
Since $D(A_\max)=\bH_{P_A^\delta}^{\mu,\delta}(X)$ and $D(A_\min)=
\bH_{P_A^{\delta+\mu-0}}^{\mu,\delta}(X)$, cf.~Theorem~\ref{grt}, and
similarly for $D(A_\max^\ast)$, $D(A_\min^\ast)$ with $P_A^\delta$ replaced
with $P_{A^\ast}^\delta$, we have to compute the induced sesquilinear form
\begin{equation}\label{seqii}
  [\;,\,]_A \colon L_{\bS^\mu}^\delta\big/L_{\bS^{\mu}}^{\delta+\mu-0}
  \times L_{\bR^\mu}^\delta\big/L_{\bR^{\mu}}^{\delta+\mu-0} \to \C.
\end{equation}
Here $\bS^\mu=\{\bs^{\mu-j}(z);\,j\in\N_0\}$ and $\bR^\mu=\{\br^{\mu-j}(z);\,
j\in\N_0\}$, where $\bs^{\mu-j}(z)= \sigma_c^{\mu-j}(A)(z)$ and
$\br^{\mu-j}(z)= \sigma_c^{\mu-j}(A^\ast)(z)$, respectively. For the relation
between $\bs^{\mu-j}(z)$, $\br^{\mu-j}(z)$, see \eqref{kily}. We will 
evaluate the sesquilinear form \eqref{seqii} using formula \eqref{gm22}.

From \eqref{seqii}, it is seen that the spaces $D(A_\max)/D(A_\min)$,
$D(A_\max^\ast)/D(A_\min^\ast)$ are invariant under the action of the operator
$T$ from \eqref{tTt}. (This result is implicitly contained in
Theorem~\ref{proper0}.)

\textbf{Step 2.}\enspace It suffices to prove \eqref{g-f} for an arbitrary
characteristic basis $\Phi_1, \dots,\Phi_e$ of the quotient
$L_{\bS^\mu}^\delta\big/ L_{\bS^{\mu}}^{\delta+\mu-0}$. For then
non-degeneracy of the sesquilinear form \eqref{seqii} and also property
\eqref{skew} follow, where the latter holds for all $\Phi\in
L_{\bS^\mu}^\delta\big/L_{\bS^{\mu}}^{\delta+\mu-0}$, $\Psi\in
L_{\bR^\mu}^\delta\big/L_{\bR^{\mu}}^{\delta+\mu-0}$. If $\Phi'_1,
\dots,\Phi'_e$ is another characteristic basis of the quotient
$L_{\bS^\mu}^\delta\big/ L_{\bS^{\mu}}^{\delta+\mu-0}$, then we have (after
renumbering if necessary)
\[
  \Phi_i' = C\,\Phi_i, \quad 1\leq i\leq e
\]
for some linear invertible operator $C\colon L_{\bS^\mu}^\delta\big/
L_{\bS^{\mu}}^{\delta+\mu-0} \to L_{\bS^\mu}^\delta\big/
L_{\bS^{\mu}}^{\delta+\mu-0}$ that commutes with $T$. Denoting by
$C^\ast$ the adjoint to $C$ with respect to the non-degenerate
sesquilinear form \eqref{seqii} ($C^\ast$ also commutes with $T$), the
conjugate characteristic basis $\Psi'_1, \dots,\Psi'_e$ to
$\Phi'_1,\dots,\Phi'_e$ is given by
\[
  \Psi_i' = (C^\ast)^{-1}\,\Psi_i, \quad 1\leq i\leq e,
\]
where $\Psi_1, \dots,\Psi_e$ is the conjugate characteristic basis to
$\Phi_1, \dots,\Phi_e$.
 
\textbf{Step 3.}\enspace Let $\Phi_1,\dots,\Phi_e$ be a characteristic basis
of $L_{\bS^\mu}^\delta \big/L_{\bS^{\mu}}^{\delta+\mu-0}$, of characteristic
$(m_1,\dots,m_e)$ say, and let $\Psi_1,\dots,\Psi_e$ be the conjugate
characteristic basis of $L_{\bR^\mu}^\delta \big/L_{\bR^{\mu}}^{\delta+\mu-0}$
according to Theorem~\ref{main1} meaning that there are corresponding versions
of Theorems~\ref{main1}, \ref{main2} valid for the Mellin symbols
$\bt^{-\mu-j}(z)$ for $0\leq j<\mu$ in the strip $\ds-\mu+j<\Re z <\ds-\mu$,
where now in \eqref{fundamental}, \eqref{graff} elements of the quotients
$L_{\bS^\mu}^\delta \big/L_{\bS^{\mu}}^{\delta+\mu-0}$, $L_{\bR^\mu}^\delta
\big/L_{\bR^{\mu}}^{\delta+\mu-0}$ enter. Likewise, we may assume that
$\Phi_1,\dots,\Phi_e$ stem (by projection) from a characteristic basis of
$L_{\bS^\mu}^\delta$ that can be extended to a complete characteristic basis
of $L_{\bS^\mu}$.

We will make this latter assumption to keep the notation from
Theorems~\ref{main1}, \ref{main2}.

For $\Phi\in L_{\bS^\mu}^\delta\big/L_{\bS^{\mu}}^{\delta+\mu-0}$,
$\Psi\in L_{\bR^\mu}^\delta\big/L_{\bR^{\mu}}^{\delta+\mu-0}$, we rewrite
\eqref{gm22} as
\begin{multline}\label{hick}
  [\Phi,\Psi]_A = -\,\sum_{k=0}^{\mu-1} \sum_{\ds-\mu+k<\Re p<\ds} \\*
  \times\, \Res_{z=p}
  \Bigl\langle\bs^{\mu-k}(z)\Phi(p)[z-p],\bbI \Psi(q+k)[z-p]\Bigr\rangle,
\end{multline}
where $q=\dim X-2\delta-\bar p-\mu$.

Now choose $\Phi$ belonging to the Jordan basis
$\Phi_1,\dots,T^{m_1-1}\Phi_1,\dots,\Phi_e, \dots,T^{m_e-1}\Phi_e $ of the
quotient $L_{\bS^\mu}^\delta\big/ L_{\bS^{\mu}}^{\delta+\mu-0}$ and $\Psi$
belonging to the conjugate Jordan basis $\Psi_1,\dots, T^{m_1-1}\Psi_1,\dots,$
$\Psi_e,\dots,T^{m_e-1}\Psi_e$ of the quotient
$L_{\bR^\mu}^\delta\big/L_{\bR^{\mu}}^{\delta+\mu-0}$. That means that
\[
  \Phi= T^i T^{m_h^{p+l+1}}\Phi_h
\]
for some $h,\,p,\,l,\,i$, where $\ds-\mu<\Re p \leq \ds-(\mu-1)$,
$\Re p+l\leq \ds$, and $0\leq i<m_h^{p+l} -m_h^{p+l+1}$. We may further
assume that
\[
  \Psi = T^jT^{m_{h^\ast}^{q+1}}\Psi_{h^\ast}, 
\]
where $q=\dim X-2\delta-\bar p-\mu$ and $0\leq j< m_{h^\ast}^{q+l+1}-
m_{h^\ast}^{q+1}$, since otherwise 
$[\Phi,\Psi]_A=0$.

Under these hypotheses, in \eqref{hick} there are non-zero residues at
most at $z=p+r$ for $r=0,\dots,l$, i.e.,
{\allowdisplaybreaks
\begin{align*}
  [\Phi,\Psi]_A &= (-1)^{j+1}\sum_{k=0}^l \sum_{r=k}^l
  \Res_{z=p+r} \left\langle\bs^{\mu-k}(z) T^i
  T^{m_h^{p+l+1}}\Phi_h(p+r)[z-p-r], \right. \\ 
  & \hspace{140pt} \left. T^j\bbI 
  T^{m_{h^\ast}^{q+1}}\Psi_{h^\ast}(q-r+k)[z-p-r]\right\rangle \\
  &= (-1)^{j+1}\sum_{k=0}^l \sum_{r=k}^l
  \Res_{z=p+r-k} \\
  & \hspace{70pt} \left\langle\bs^{\mu-k}(z+k) T^i
  T^{m_h^{p+l+1}}\Phi_h(p+r)[z-p-r+k], \right. \\ 
  & \hspace{140pt} \left. T^j\bbI 
  T^{m_{h^\ast}^{q+1}}\Psi_{h^\ast}(q-r+k)[z-p-r+k]\right\rangle \\
  &= (-1)^{j+1}\sum_{k=0}^l \sum_{r=0}^{l-k}
  \Res_{z=p+r} \\
  & \hspace{70pt} \left\langle\bs^{\mu-k}(z+k) T^i
  T^{m_h^{p+l+1}}\Phi_h(p+r+k)[z-p-r], \right. \\ 
  & \hspace{140pt} \left. T^j\bbI 
  T^{m_{h^\ast}^{q+1}}\Psi_{h^\ast}(q-r)[z-p-r]\right\rangle \\
  &= (-1)^{j+1}\sum_{r=0}^l 
  \Res_{z=p+r} \\
  & \hspace{70pt} \left\langle\sum_{k=0}^{l-r} \bs^{\mu-k}(z+k) T^i
  T^{m_h^{p+l+1}}\Phi_h(p+r+k)[z-p-r], \right. \\ 
  & \hspace{140pt} \left. T^j\bbI 
  T^{m_{h^\ast}^{q+1}}\Psi_{h^\ast}(q-r)[z-p-r]\right\rangle \\
  &= (-1)^{j+1} \sum_{r=0}^l 
  \Res_{z=p+r} (z-p-r)^{i+j}\left\langle \bbTh_{l-r}
  (T^{m_h^{p+l+1}}\Phi_h)[z], \right. \\ 
  & \hspace{140pt} \left. \bbI 
  T^{m_{h^\ast}^{q+1}}\Psi_{h^\ast}(q-r)[z-p-r]\right\rangle \\
  &= (-1)^{j+1} \Res_{z=p} 
  (z-p)^{i+j} \sum_{r=0}^l \left\langle \bbTh_{l-r}
  (T^{m_h^{p+l+1}}\Phi_h)[z+r], \right. \\ 
  & \hspace{140pt} \left. \bbI 
  T^{m_{h^\ast}^{q+1}}\Psi_{h^\ast}(q-r)[z-p]\right\rangle \\
\end{align*}}
Therefore, 
\[
  [\Phi,\Psi]_A = \begin{cases}
  (-1)^{j+1} & \text{if $\tau^\ast(h)=h^\ast$, $i+j=m_h^p-m_h^{p+l+1}-1$,} \\
  0 & \text{otherwise,}
  \end{cases}
\]
by virtue of Theorem~\ref{main2}.

This completes the proof.
\end{proof}


\section{Examples}\label{sec7}


We discuss two examples of ordinary differential operators on the
half-line~$\Rp$. The first example demonstrates the usage of
Theorem~\ref{main} for the computation of the boundary sequilinear form, while
in the second example it is shown how our fundamental formulas like
\eqref{fundamental} can be independently verified.


\subsection{First example}

This example concerns the cone-degenerate third-order operator
\[
  A = \partial_t^3 + t^{-1}\partial_t^2 \quad \text{on $\Rp$.}
\]
The conormal symbols are 
\[
  \sigma_c^3(A)(z)=-z(z+1)^2
\] 
and $\sigma_c^{3-j}(A)(z)=0$ for $j\geq1$. Thus, $1$, $t\log t$, and $t$ are
exact solutions to the equation $Au=0$. A complete characteristic basis
$\Phi_1,\Phi_2$ of $L_A$ is given by
\[
  \Phi_1(0)=(1), \quad \Phi_2(-1)=(1,0), 
\]
and $\Phi_1(p)=0$ for $p\neq0$, $\Phi_2(p)=0$ for $p\neq-1$. 

We choose \underline{$\delta=-1$}. Then we have $A^\ast= -\partial_t^3 -
5t^{-1}\partial_t^2 - 4t^{-2}\partial_t$, $\sigma_c^3(A^\ast)(z)=z(z-1)^2$,
and $\sigma_c^{3-j}(A^\ast)(z)=0$ for $j\geq1$. From
\[
  -\,\frac{1}{z(z+1)^2} = \frac{1}{(z+1)^2} + \frac{1}{z+1} - \frac{1}{z} 
\]
we infer that the complete characteristic basis $\Psi_1,\Psi_2$ of
$L_{A^\ast}$ that is conjugate to $\Phi_1,\Phi_2$ is given by
\[
  \Psi_1(1)=(1,-1), \quad \Psi_2(0)=(1), 
\]
and $\Psi_1(p)=0$ for $p\neq1$, $\Psi_2(p)=0$ for $p\neq0$, where
$\tau^\ast(1)=2$, $\tau^\ast(2)=1$.

Writing 
\begin{align*}
  u(t) &=\omega(t)\left(\alpha + \beta_0 t\log t + \beta_1 t \right) 
  + u_0(t), \\
  v(t) &=\omega(t)\left(\gamma_0 t^{-1}\log t + \gamma_1 t^{-1} +
  \delta\right) + v_0(t)
\end{align*}
for $\alpha,\,\beta_0,\,\beta_1,\,\gamma_0,\,\gamma_1,\,\delta\in\C$, where
$\omega(t)$ is a cut-off function and $u_0\in D(A_\min)$, $v_0\in
D(A_\min^\ast)$, we then obtain
\[
  [u,v]_A = - \alpha\bar \delta + \beta_0\bar\gamma_0 +
  \beta_0 \bar\gamma_1 - \beta_1 \bar\gamma_0.
\]


\subsection{Second example}

We consider the non-degenerate, second-order, constant coefficient
operator 
\[
  A = \partial_t^2 + a\partial_t + b \quad \text{on $\Rp$,}
\]
where $a,\,b\in\C$. We have $A^\ast=\partial_t^2-\bar a\partial_t+\bar b$
(with \underline{$\delta=0$}) and the Green's formula is directly checked to
be
\begin{equation}\label{exam1}
  [u,v]_A = u(0)\bar v'(0) - u'(0) \bar v(0) - a\, u(0) \bar v(0),
  \quad u,\,v\in {\mathcal S}(\oRp).
\end{equation}

The space $L_A^0$ has characteristic $\{(-j; 1,1,\dots);\,j=0,1\}$.
Therefore, $L_A^{2-0}=\{0\}$, the quotient $D(A_\max)/$ $D(A_\min) \cong
L_A^0$ is two-dimensional, and the elements of $D(A_\max)/D(A_\min)$ are in
one-to-one correspondence with the (analytic, as turns out) solutions
$u(t)=u(t;\alpha,\beta)$ for $\alpha,\,\beta\in\C$ to
\[
   Au=0, \quad u(0)=\alpha, \quad u'(0)=\beta.
\]
In the following, we shall assume this identification to be made.

A complete characteristic basis of $ L_A= L_A^0$ is given by
\begin{align*}
  u_1(t)=u(t;1,0), \quad & u_2(t)=u(t;0,1). \\
\intertext{We look at \eqref{exam1} to find the conjugate
    complete characteristic basis of $ L_{A^\ast}$ to be}
  v_1(t)=u(t;1,\bar a), \quad & v_2(t)=v(t;0,-1),
\end{align*}
where $v(t)=v(t;\alpha,\beta)$ is the solution to $A^\ast v=0$,
$v(0)=\alpha$, $v'(0)=\beta$.

\begin{proposition}
\textup{(a)} We have 
\begin{align*}
  & u_1(t) = 1 + \sum_{j\geq 2} (-1)^{j-1}\,\frac{b\Pi_{j-2}(a,b)}{j!}\,t^j, 
  \quad 
  u_2(t) = \sum_{j\geq 1} (-1)^{j-1}\,\frac{\Pi_{j-1}(a,b)}{j!}\,t^j \\
\intertext{and}
  & v_1(t) = \sum_{j\geq 0} \frac{\Pi_j(\bar a,\bar b)}{j!}\,t^j, \quad 
  v_2(t) = - \sum_{j\geq 1} \frac{\Pi_{j-1}(\bar a,\bar b)}{j!}
  \,t^j,
\end{align*}
where $\Pi_0(a,b)=1$, $\Pi_1(a,b)=a$, 
\[
  \Pi_j(a,b) = a\Pi_{j-1}(a,b) - b\Pi_{j-2}(a,b), \quad j=2,3,\dots
\]
\textup{(}i.e., $\Pi_2(a,b) = a^2-b$, $\Pi_3(a,b)=a^3-2ab$,
$\Pi_4(a,b)=a^4-3a^2b+b^2$, etc.\/\textup{)}.

\textup{(b)} We also have 
\[
  \bt^{-k-2}(z+2) = \frac{\Pi_k(a,b)}{(z-k)(z-k+1)\dots z(z+1)},
  \quad k=0,1,2,\dots,
\]
where $\bt^{-k-2}(z)$ has the same meaning as before. In particular, the poles
of $\bt^{-k-2}(z)$ are simple and, for $l=-1,0,1,\dots,k$,
\[
  \operatorname{Res}_{z=l} \bt^{-k-2}(z+2) = 
  \frac{(-1)^{k-l}\Pi_k(a,b)}{(k-l)!(l+1)!}.
\]
\end{proposition}

The key in re-proving formulas like \eqref{fundamental} is:

\begin{lemma}
\textup{(i)} For $l\geq2$, $0\leq j\leq l-2$, 
\[
  \Pi_l(a,b) = \Pi_{j+1}(a,b)\Pi_{l-j-1}(a,b) - b\Pi_j(a,b) \Pi_{l-j-2}(a,b).
\]

\textup{(ii)} For $j\geq0$, 
\[ 
  \Pi_j(-\bar a,\bar b)=(-1)^j\,\overline{\Pi_j(a,b)}.
\]
\end{lemma}


\appendix
\section{Local asymptotic types}\label{appA}


We briefly discuss the notion of local asymptotic type, i.e., asymptotic types
at one singular exponent $p\in\C$ in \eqref{asymp}. Moreover, we investigate
an analogue of the boundary sesquilinear form in this simpler situation, see
\eqref{harry}. Most of the material is taken from \textsc{Witt} \cite{Wit02}.

Let $E$ be a Banach space, $E'$ be its topological dual, and let
$\langle\;,\,\rangle$ denote the dual pairing between $E,\,E'$. Pick $p\in\C$.
For notations like $\cM_p(E)$, $\cA_p(E)$, $E^\infty$, the right-shift operator
$T$ acting on $E^\infty$, $\Phi\otimes\Psi[z-p]$ for $\Phi\in E^\infty$, 
$\Psi\in E'{}^\infty$, and the identification $\cM_p(E)\big/\cA_p(E)\cong
E^\infty$, see Section~\ref{nota}. 

Let $\cM_p^\fin(\cL(E))$ be the space of germs of $\cL(E)$--valued 
finitely meromorphic functions $F(z)$ at $z=p$, i.e.,
\begin{equation}\label{germ}
  F(z) = \frac{F_0}{(z-p)^\nu} + \frac{F_1}{(z-p)^{\nu-1}} + \dots +
  \frac{F_{\nu-1}}{z-p} + \sum_{j\geq 0} F_j (z-p)^j,
\end{equation}
where $F_0,F_1,\dots,F_{\nu-1}\in \cL(E)$ are finite-rank operators. Let
$\cM_p^\nor(\cL(E))$ be the space of germs of $\cL(E)$-valued normally
meromorphic functions $F(z)$ at $z=p$, i.e., the space of finitely meromorphic
functions $F(z)$, where, in addition, $F(z)$ for $z\neq p$ close to $p$ is
invertible and $F_\nu\in\cL(E)$ is a Fredholm operator. $\cM_p^\nor(\cL(E))$
is the group of invertible elements of the algebra $\cM_p^\fin(\cL(E))$.

For $F\in\cA_p(\cL(E))$, let $L_F$ denote the space of all $(\phi_0,\phi_1,
\dots,\phi_{m-1})\in E^\infty$ such that
\[
  F(z) \left(\frac{\phi_0}{(z-p)^m} + \frac{\phi_1}{(z-p)^{m-1}} + \dots +
  \frac{\phi_{m-1}}{(z-p)}\right) \in \cA_p(E).
\]

\begin{remark}
The theory can be also developed for $F\in\cM_p^\fin(\cL(E))$ upon an
appropriate modification of the definition of $L_F$. For $F\in
\cM_p^\nor(\cL(E))$, $L_F$ is again an asymptotic type, and Propositions
\ref{awq}, \ref{awq2}, and \ref{awq3} continue to hold in this case. See
\textsc{Witt}~\cite{Wit02}.
\end{remark}

\begin{definition}
A local asymptotic type $J\subset E^\infty$ is a finite-dimensional linear
subspace that is invariant under the action of the right shift operator $T$.
The set of all local asymptotic types is denoted by $\cJ(E)$.
\end{definition}

Note that $T$ as acting on $J$ is nilpotent. The characteristic
$(m_1,\dots,m_e)$ of $T$ on $J$ is called the characteristic of the asymptotic
type.

\begin{proposition}
We have
\[
  \cJ(E) = \bigl\{L_F\bigm|F\in \cA_p(\cL(E))\cap \cM_p^\nor(\cL(E))\bigr\}.
\]
\end{proposition}

\begin{proposition}\label{awq}
For $F\in \cA_p(\cL(E))\cap \cM_p^\nor(\cL(E))$, we have $F^t\in \cA_p(\cL(E))
\cap \cM_p^\nor(\cL(E'))$. Moreover, for each characteristic basis
$\Phi_1,\dots,\Phi_e$ of $L_F$, there exists a unique characteristic basis
$\Psi_1,\dots,\Psi_e$ of $L_{F^t}$ such that
\[
  \bigl[F^{-1}(z)\bigr]_p^\ast = \sum_{i=1}^e (\Phi_i\otimes\Psi_i)[z-p].
\]
In particular, both asymptotic types $L_F$, $L_{F^t}$ have the same
characteristic. 
\end{proposition}

The next result is Keldysh's formula, cf.~\textsc{Keldysh} \cite{Kel51},
\textsc{Kozlov--Maz'ya} \cite{KM99}.

\begin{proposition}\label{awq2}
For $\Phi_1,\dots,\Phi_e$ and $\Psi_1,\dots,\Psi_e$ as in\/ 
\textup{Proposition~\ref{awq}},
\begin{equation}\label{perw}
  \bigl\langle F(z)\Phi_i[z-p],\Psi_j[z-p]\bigr\rangle = 
  \delta_{ij}(z-p)^{-m_i} + O(1) \quad \textnormal{as $z\to p$.}
\end{equation}
\end{proposition}

\begin{remark}
Writing $F(z)$ as in \eqref{germ}, \eqref{perw} can be rewritten as 
\[
  \sum_{\nu+r+s=m_i+l}
  \bigl\langle F_\nu\phi_r^{(i)},\psi_s^{(j)}\bigr\rangle = 
  \delta_{ij}\delta_{0l} 
\] 
for $0\leq l\leq m_j-1$, where summation is restricted to the range $0\leq
r\leq m_i-1$, $0\leq s\leq m_j-1$, and
$\Phi_i=\bigl(\phi_0^{(i)},\phi_1^{(i)}, \dots,\phi_{m_i-1}^{(i)}\bigr)$,
$\Psi_j=\bigl(\psi_0^{(j)},\psi_1^{(j)}, \dots,\psi_{m_j-1}^{(j)}\bigr)$.
\end{remark}

For $F\in \cA_p(\cL(E))\cap \cM_p^\nor(\cL(E))$, we finally introduce the
bilinear form $[\;,\,]_F$ defined on the product $L_F\times L_{F^t}$ by
\begin{equation}\label{harry}
  [\Phi,\Psi]_F := \Res_{z=p}\,\bigl\langle F(z)\Phi[z-p],\Psi[z-p]
  \bigr\rangle.
\end{equation}

\begin{proposition}\label{awq3}
Evaluated on the bases $T^r\Phi_i$ for $1\leq i\leq e$, $0\leq r\leq m_i-1$ of
$L_F$ and $T^s\Psi_j$ for $1\leq j\leq e$, $0\leq s\leq m_j-1$ of $L_{F^t}$,
\[
  \bigl[T^r\Phi_i,T^s\Psi_j\bigr]_F = \begin{cases}
  1 & \textnormal{if $i=j$, $r+s=m_i-1$,} \\
  0 & \textnormal{otherwise.}
  \end{cases}
\]
\end{proposition}
\begin{proof}
This follows immediately from Proposition~\ref{awq2}.
\end{proof}


\section{Function spaces with asymptotics}\label{appB}


The maximal and minimal domains of cone-degenerate elliptic differential
operators are cone Sobolev spaces with asymptotics, as we are going to
demonstrate now. We refer to \textsc{Schulze}~\cite{Sch91, Sch98} for more on
function spaces with asymptotics, where, however, asymptotics are observed on
so-called ``half-open weight intervals,'' a setting leading to Fr\'echet
spaces. The present setting due to \textsc{Liu--Witt} \cite{LW01}, where
asymptotics are observed on ``closed weight intervals,'' provides a scale of
Hilbert spaces.


\subsection{Weighted cone Sobolev spaces}

Let $Mu(z)=\tilde u(z)=\int_0^\infty t^{z-1}u(t)\,dt$ for $z\in\C$ (or subsets
thereof) be the Mellin transformation, suitably extended to certain
distribution classes. Recall that
\begin{equation}\label{mttt}
  M\colon L^2(\Rp,t^{-2\delta}dt) \to L^2(\Gamma_{1/2-\delta};
  (2\pi i)^{-1}dz),
\end{equation}
is an isometry, where $\Gamma_\gamma:= \{z\in\C\,|\,\Re z=\gamma\}$ for
$\gamma\in\R$. Moreover,
\begin{align*}
  & M_{t\to z}\bigl\{(-t\partial_t)u\bigr\}(z) = z\,\tilde u(z), \\
  & M_{t\to z}\bigl\{t^{-p}u\bigr\}(z) = \tilde u(z-p), \quad p\in\C.
\end{align*}
The function 
\[
  \bm_{p,k}(z,y) := M_{t\to z}\left\{\frac{(-1)^k}{k!}\,\omega(t) 
  t^{-p} \log^k t\,\phi(y)\right\}, 
\]
where $p\in\C$, $k\in\N_0$, $\phi\in C^\infty(Y)$, and $\omega(t)$ is a
cut-off function, belongs to $\cM_\as^{-\infty}(Y)$. Furthermore,
\[
  \bm_{p,k}(z)-\frac{\phi(y)}{(z-p)^{k+1}}\in\cA(\C;C^\infty(Y)).
\]

For $s,\,\delta\in\R$, the Hilbert space $\cH^{s,\delta}(X)$ consists of all
$u\in H_{\loc}^s(X^\circ)$ such that $M_{t\to z}\{\omega u\}(z)\in
L_{\loc}^2\bigl( \Gamma_\ds;H^s(Y)\bigr)$ and
\[
   \frac{1}{2\pi i}\int_{\Gamma_\ds}\bigl\|R^s(z)
   M_{t\to z}\{\omega u\}(z)\bigr\|_{L^2(Y)}^2\,dz < \infty.
\]
Here, $R^s(z)\in L^s(Y;\Gamma_\ds)$ is an order-reducing family, i.e.,
$R^s(z)$ is para\-meter-dependent elliptic and $R^s(z)\colon H^{s+s'}(Y)\to
H^{s'}(Y)$ is invertible for all $s'\in\R$, $z\in \Gamma_\ds$. For instance,
if $\bm(z)\in\cM_{\as}^s(Y)$ is elliptic and the line $\Gamma_\ds$ is free of
poles of $\bm(z)$, then $\bm(z)\big|_{\Gamma_\ds}$ is such an order-reduction.


\subsection{Cone Sobolev spaces with asymptotics}

The starting point is the following observation:

\begin{theorem}[\textsc{Liu--Witt} {\cite[Theorem~2.43]{LW01}}]\label{oben}
Let $s,\,\delta\in\R$ and $P\in\As^\delta(Y)$ be a proper asymptotic type.
Then there exists an elliptic Mellin symbol $\bm_P^s(z)\in\cM_\cO^s(Y)$ such
that the line $\Gamma_\ds$ is free of poles of $\bm_P^s(z)^{-1}$ and, for
$\mathfrak{S}^s= \bigl\{\bm_P^s(z),0,0,\dots \bigr\}\in \Symb_M^s(Y)$,
$L_{\mathfrak{S}^s}^\delta$ represents the asymptotic type $P$.
\end{theorem}

\begin{definition}\label{ght}
Let $s\geq0$, $\delta\in\R$, and $P\in\As^\delta(Y)$ be proper. Then the space
$\bH_P^{s,\delta}(X)$ consists of all functions $u\in \cH^{s,\delta}(X)$ such
that $M_{t\to z}\{\omega u\}(z)$ is meromorphic for $\Re z>\ds-s$ with values
in $H^s(Y)$,
\[
  \bm_P^s(z)M_{t\to z}\{\omega u\}(z)\in \cA\bigl(
  \{z\in\C\,|\,\Re z>\ds-s\};L^2(Y)\bigr),
\]
where $\bm_P^s(z)$ is as in \textnormal{Theorem~\ref{oben}}, and 
\[
  \sup_{0<s'<s} \frac{1}{2\pi i}
  \int_{\Gamma_{\ds-s'}}\bigl\|\bm_P^s(z)M_{t\to z}\{\omega u\}(z)
  \bigr\|_{L^2(Y)}^2\,dz < \infty.
\]
\end{definition}

We list some properties of the spaces $\bH_P^{s,\delta}(X)$:

\begin{proposition}\label{prop_csa}
\textup{(a)} $\bigl\{\bH_P^{s,\delta}(X);\,s\geq 0\bigr\}$ is an interpolation
scale of Hilbert spaces with respect to the complex interpolation method.

\textup{(b)} $\bH_\cO^{s,\delta}(X)= \cH^{s,\delta+s}(X)$.

\textup{(c)} We have
\begin{multline*}
  \bH_P^{s,\delta}(X) = \bH_\cO^{s,\delta}(X) 
  \oplus\Bigl\{\omega(t)\,\sum_{
  \Re p>\ds-s}
  \sum_{k+l=m_p-1}\frac{(-1)^k}{k!}\,t^{-p}\log^kt\,\phi_l^{(p)}(y) \Bigm| \\
  \text{$\Phi(p)=(\phi_0^{(p)},\dots,\phi_{m_p-1}^{(p)})$ for some 
    $\Phi\in J$}\Bigr\},
\end{multline*}
where the linear space $J\subset\cE_V^\delta(Y)$ represents the asymptotic
type $P$, provided that
\[
  \Re p\neq\ds-s, \quad p\in V.
\]

\textup{(d)} $\bH_P^{s,\delta}(X)\subseteq \bH_{P'}^{s',\delta'}(X)$ if and
only if $s\geq s'$, $\delta+s\geq\delta'+s'$, and $P\preccurlyeq P'$ up to the
conormal order $\delta'+s'$.

\textup{(e)} $C_P^\infty(X) := \bigcap_{s\geq0} \bH_P^{s,\delta}(X)$ is dense
in $\bH_P^{s,\delta}(X)$.
\end{proposition}

\begin{proposition}\label{explai}
The spaces $\bH_P^{s,\delta}(X)$ are invariant under coordinate changes in the
sense explained in\/ \textup{Remark~\ref{ccc}}.
\end{proposition}


\subsection{Mapping properties and elliptic regularity}

Here we are concerned with the regularity and asymptotics of solutions $u$ to
the equation
\begin{equation}\label{hrt}
  Au(x) = f(x) \quad \text{on $X^\circ$},
\end{equation}
where $A\in \DiffF^\mu(X)$ is elliptic. Assuming $u\in \bH^{0,\delta}(X)$ and
$f\in \bH_Q^{s,\delta}(X)$, where $s\geq0$ and $Q\in\As^\delta(Y)$, we are
going to show that $u\in \bH_P^{s+\mu,\delta}(X)$ for some resulting
$P\in\As^\delta(Y)$. By interior elliptic regularity, we already know that
$u\in H_\loc^{s+\mu}(X^\circ)$. So we are left with the behavior of $u=u(x)$
as $x\to\partial X$.

Let $P_A^\delta$ be the asymptotic type represented by $L_A^\delta$.
Similarly, let $P_A^{\delta+\mu-0}\preccurlyeq P_A^\delta$ be the asymptotic
type represented by $L_A^{\delta+\mu-0}$. Then $P_A^{\delta+\mu-0}$ is the
largest asymptotic type that coincides with the empty asymptotic type, $\cO$,
up to the conormal order $\delta+\mu-0$. Note that, for each
$P\in\As^\delta(Y)$ satisfying $P\preccurlyeq P_A^\delta$ up to the conormal
order $\delta+\mu$, there is a $Q\in\As^\delta(Y)$ such that
\[
  A\colon \bH_P^{s+\mu,\delta}(X) \to \bH_Q^{s,\delta}(X)
\]
for all $s\geq0$. The minimal such $Q\in\As^\delta(Y)$ is denoted by
$Q^\delta(P;A)$. In particular, $Q^\delta(P_A;A)=Q^\delta(\cO;A)=\cO$.

The question raised for equation \eqref{hrt} is answered by the next result:

\begin{proposition}
Let $A\in \DiffF^\mu(X)$ be elliptic. Then\/\textup{:} 

\textup{(a)} The map
\begin{multline}\label{grt0}
  \bigl\{P\in\As^\delta(Y)\colon P\succcurlyeq P_A^\delta,\,
  \text{$P$ coincides with $P_A^\delta$ up to} \\
  \text{the conormal order $\delta+\mu$}\bigr\} \to \As^\delta(Y),
  \quad P\mapsto Q^\delta(P;A)
\end{multline}
is an order-preserving bijection.

\textup{(b)} For any solution $u$ to \eqref{hrt}, $u\in
\bH^{0,\delta}(X)$ and $f\in \bH_Q^{s,\delta}(X)$ implies $u\in
\bH_{P^\delta(Q;A)}^{s+\mu,\delta}(X)$, where $Q\mapsto P^\delta(Q;A)$ is the
inverse to \eqref{grt0}.
\end{proposition}

Note that $Q^\delta(P^\delta(Q;A);A)=Q$. Therefore, $P\mapsto P^\delta(
Q^\delta(P;A);A)$ is a hull operation. Note also that both maps $P\to
Q^\delta(P;A)$ and $Q\mapsto P^\delta(Q;A)$ to \eqref{grt0} can be computed
purely on the level of the complete conormal symbols
$\{\sigma_c^{\mu-j}(A)(z);\,j\in\N_0\}$.

\begin{theorem}\label{grt}
Let $A\in \DiffF^\mu(X)$ be elliptic. Then
\begin{equation}\label{grt1}
  D(A_\max) = \bH_{P_A^\delta}^{\mu,\delta}(X), \quad D(A_\min) = 
  \bH_{P_A^{\delta+\mu-0}}^{\mu,\delta}(X).
\end{equation}
In particular, 
\begin{equation}\label{grt2}
  D(A_\max)\big/D(A_\min) \cong L_A^\delta\big/L_A^{\delta+\mu-0}.
\end{equation}
\end{theorem}
\begin{proof}
\eqref{grt1} is a consequence of elliptic regularity, while \eqref{grt2}
follows from the description given in Proposition~\ref{prop_csa} (c) 
and interpolation. 
\end{proof}


\nocite{*}
\bibliographystyle{amsplain}
\bibliography{green}

\end{document}